\documentclass{gtart_h}  

\def\ifplaintex{\expandafter\ifx\csname documentclass\endcsname\relax}

\def\ifplaintex{\expandafter\ifx\csname documentclass\endcsname\relax}

\ifplaintex 
\else
\fi

\expandafter\ifx\csname epsfbox\endcsname\relax\input epsf\fi

\def\gt{{\mathsurround=0pt\it $\cal G\mskip-2mu$eometry \&\ 
$\cal T\!\!$opology}}        

\def\gtp{{\mathsurround=0pt\it $\cal G\mskip-2mu$eometry \&\ 
$\cal T\!\!$opology $\cal P\!$ublications}}  

\def\lognumber#1{\def\thelognumber{#1}}
\def\volumenumber#1{\def\thevolumenumber{#1}}
\def\papernumber#1{\def\thepapernumber{#1}}
\def\volumeyear#1{\def\thevolumeyear{#1}}

\def\pagenumbers#1#2{\def\startpage{#1}\def\finishpage{#2}}
\def\published#1{\def\publishdate{#1}}
\def\proposed#1{\def\theproposer{#1}}
\def\seconded#1{\def\theseconders{#1}}
\def\received#1{\def\receiveddate{#1}}

\def\accepted#1{\def\accepteddate{#1}}

\def\asciiaddress#1{\def\theasciiaddress{#1}}
\def\asciiemail#1{\def\theasciiemail{#1}}

\long\def\asciiabstract#1{\long\def\theasciiabstract{#1}}

\let\\\par\let\thelognumber\relax
\let\thevolumenumber\relax\let\thepapernumber\relax
\let\thevolumeyear\relax\let\thesamplenumber\relax\let\startpage\relax
\let\finishpage\relax\let\publishdate\relax\let\receiveddate\relax
\let\reviseddate\relax\let\accepteddate\relax\let\theasciititle\relax
\let\theasciiauthors\relax\let\theasciiaddress\relax
\let\theasciiabstract\relax
\let\theasciiemail\relax\let\theshortauthors\relax\let\theshorttitle\relax

\long\def\maketitlep{   

\count0=\startpage

\gt\hfill      
\hbox to 77pt{\vbox to 0pt{\vglue -15pt\epsfbox{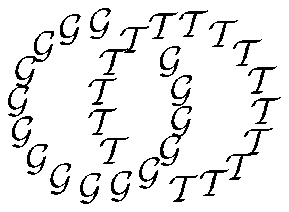}\vss}\hss}
\break
{\small\ifx\thesamplenumber\relax 
Volume \else Sample
\fi\thevolumenumber\ (\thevolumeyear)
\startpage--\finishpage\nl
Published: \publishdate}
\vglue 0.5truein plus 0.4fil minus 0.1truein

{\parskip=0pt\leftskip 0pt plus 1fil\def\\{\par\smallskip}{\ifplaintex\large
\else\Large\fi\bf\thetitle}\par\medskip}   

\vglue 0pt plus 0.1fil 

{\parskip=0pt\leftskip 0pt plus 1fil\def\\{\par}{\sc\theauthors}
\par\medskip}

\vglue 0pt plus 0.1fil 

{\small\parskip=0pt\let\newline\\
{\leftskip 0pt plus 1fil\def\\{\par}{\sl\theaddress}\par}
\expandafter\ifx\theemail\relax    
\relax\else\vglue 5pt plus 0.02fil minus 2pt\def\\{\stdspace{\rm 
and}\stdspace} 
\cl{Email:\stdspace\tt\theemail}\fi
\ifx\theurl\relax                  
\relax\else\vglue 5pt plus 0.02fil minus 2pt\def\\{\stdspace{\rm 
and}\stdspace}
\cl{URL:\stdspace\tt\theurl}\fi\par}

\vglue 7pt plus 0.3fil minus 3pt

{\bf Abstract}
\vglue 5pt plus 0.1fil minus 2pt

\theabstract

\vglue 7pt plus 0.3fil minus 3pt

{\bf AMS Classification numbers}\quad Primary:\quad \theprimaryclass

Secondary:\quad \thesecondaryclass

\vglue 5pt plus 0.3fil minus 2pt

{\bf Keywords:}\quad \thekeywords

\vglue 10pt plus 0.5fil minus 5pt

{\small  Proposed: \theproposer\hfill Received: \receiveddate\nl
Seconded: \theseconders\hfill 
\ifx\reviseddate\relax                         
Accepted: \accepteddate                        
\else
Revised: \reviseddate                          
\fi}
\eject
}       

\font\phead=cmsl9 scaled 950
\font=cmsl9 scaled 1050
\font\pnum=cmbx10 scaled 913
\font\lnum=cmbx10 
\font\pfoot=cmsl9 scaled 950
\font=cmsl9 scaled 1050
\ifplaintex
\headline{\vbox to 0pt{\vskip -4.5mm\line{\small\phead\ifnum
\count0=\startpage ISSN 1364-0380 (on line)
1465-3060 (printed) \hfill {\pnum\folio}\else\ifodd\count0\def\\{ }%
\ifx\theshorttitle\relax\thetitle\else\theshorttitle\fi\hfill{\pnum\folio}
\else\def\\{ and }{\pnum\folio}\hfill\ifx\theshortauthors\relax\theauthors
\else\theshortauthors\fi\fi\fi}\vss}}
\footline{\vbox to 0pt{\vglue 0mm\line{\small\pfoot\ifnum\count0=\startpage
\copyright\ \gtp\hfill\else
\gt, Volume \thevolumenumber\ (\thevolumeyear)\hfill\fi}\vss
}}
\else
\makeatletter
\def\@oddhead{{\small\lhead\ifnum\count0=\startpage ISSN 1364-0380 (on line)
1465-3060 (printed) \hfill {\lnum\number\count0}\else\ifodd\count0
\def\\{ }\ifx\theshorttitle\relax \thetitle \else\theshorttitle\fi\hfill
{\lnum\number\count0}\else\def\\{ and }{\lnum\number\count0}
\hfill\ifx\theshortauthors\relax 
\theauthors\else\theshortauthors\fi\fi\fi}}\def\@evenhead{\@oddhead}
\def\@oddfoot{\small\lfoot\ifnum\count0=\startpage\copyright\ \gtp\hfill\else
\gt, Volume \thevolumenumber\ (\thevolumeyear)\hfill\fi}
\def\@evenfoot{\@oddfoot}
\makeatother
\fi

\newwrite\gtoutfile
\long\gdef\makeheadfile{  
{\def\\{, }\def\s{ }
\immediate\openout\gtoutfile head.xxx
\immediate\write\gtoutfile{Proxy-for: \ifx\theasciiauthors\relax
\theauthors\else\theasciiauthors\fi\s<\ifx\theasciiemail\relax\theemail\else\theasciiemail\fi>}
\immediate\write\gtoutfile{\noexpand\\}
\immediate\write\gtoutfile{Authors: \ifx\theasciiauthors\relax
\theauthors\else\theasciiauthors\fi}
{\def\\{ }\immediate\write\gtoutfile{Title: \ifx\theasciititle\relax
\thetitle\else\theasciititle\fi}}
\immediate\write\gtoutfile{Subj-class: GT or SG or MG etc}
\immediate\write\gtoutfile{MSC-class: \theprimaryclass\ifx\thesecondaryclass\relax\else, \thesecondaryclass\fi}
\immediate\write\gtoutfile{Journal-ref: Geom. Topol. \thevolumenumber
(\thevolumeyear) \startpage-\finishpage}
\immediate\write\gtoutfile{Comments: Published by Geometry and Topology at}
\immediate\write\gtoutfile{\s\s http://www.maths.warwick.ac.uk/gt/GTVol\thevolumenumber/paper\thepapernumber.abs.html}
\immediate\write\gtoutfile{\noexpand\\}
\immediate\write\gtoutfile{}
\ifx\theasciiabstract\relax
\immediate\write\gtoutfile{\theabstract}\else
\immediate\write\gtoutfile{\theasciiabstract}\fi
\immediate\write\gtoutfile{}
\immediate\write\gtoutfile{\noexpand\\}
\immediate\write\gtoutfile{}
\immediate\closeout\gtoutfile}}  

\def\maketitlepage{\maketitlep\makeheadfile}
\let\maketitle\maketitlepage

\lognumber{486}
\received{12 August 2004}
\volumenumber{9}\papernumber{25}\volumeyear{2005}
\pagenumbers{1115}{1146}   
\published{1 June 2005}
\accepted{7 May 2005}
\proposed{Rob Kirby}
\seconded{Simon Donaldson, Gang Tian}

\usepackage{amsmath, amssymb}
\usepackage[arrow,matrix]{xy}

\def\eq#1{{\rm(\ref{#1})}}
\theoremstyle{plain}
\newtheorem{thm}{Theorem}[section]
\newtheorem{lem}[thm]{Lemma}
\newtheorem{prop}[thm]{Proposition}

\theoremstyle{definition}
\newtheorem{dfn}[thm]{Definition}
\newtheorem{rem}[thm]{Remark}
\def\Ker{\mathop{\rm Ker}}
\def\Coker{\mathop{\rm Coker}}
\def\ind{\mathop{\rm ind}}
\def\Re{\mathop{\rm Re}}
\def\vol{\mathop{\rm vol}}
\def\SO{\mathbin{\rm SO}}

\def\Hol{{\textstyle\mathop{\rm Hol}}}
\def\ge{\geqslant}
\def\le{\leqslant}
\def\C{{\mathbin{\mathbb C}}}
\def\R{{\mathbin{\mathbb R}}}

\def\Z{{\mathbin{\mathbb Z}}}
\def\D{{\mathbin{\mathcal D}}}
\def\H{{\mathbin{\mathcal H}}}
\def\M{{\mathbin{\mathcal M}}}
\def\al{\alpha}
\def\be{\beta}
\def\ga{\gamma}
\def\de{\delta}
\def\ep{\epsilon}

\def\la{\lambda}
\def\ze{\zeta}
\def\th{\theta}
\def\vp{\varphi}
\def\si{\sigma}

\def\om{\omega}
\def\De{\Delta}
\def\Ga{\Gamma}
\def\Th{\Theta}
\def\La{\Lambda}
\def\Om{\Omega}
\def\ts{\textstyle}
\def\sst{\scriptscriptstyle}
\def\sm{\setminus}
\def\na{\nabla}
\def\pd{\partial}
\def\op{\oplus}
\def\ot{\otimes}
\def\bigop{\bigoplus}
\def\iy{\infty}
\def\ra{\rightarrow}
\def\longra{\longrightarrow}

\def\t{\times}
\def\w{\wedge}
\def\d{{\rm d}}

\def\ci{\circ}
\def\ti{\widetilde}

\def\bar{\overline}
\def\md#1{\vert #1 \vert}
\def\nm#1{\Vert #1 \Vert}
\def\bmd#1{\big\vert #1 \big\vert}
\def\lnm#1#2{\Vert #1 \Vert_{L^{#2}}} 
\def\bnm#1{\bigl\Vert #1 \bigr\Vert}

\begin{document}

\title[Deformations of AC coassociative 4--folds]{Deformations of
asymptotically cylindrical\\coassociative submanifolds with fixed boundary}
\authors{Dominic Joyce\\Sema Salur}
\address{Lincoln College, Oxford, OX1 3DR, UK}
\secondaddress{Department of Mathematics, Northwestern University, IL 60208, USA}

\asciiaddress{Lincoln College, Oxford, OX1 3DR, UK\\and\\Department of 
Mathematics, Northwestern University, IL 60208, USA}

\gtemail{\mailto{dominic.joyce@lincoln.oxford.ac.uk}, 
\mailto{salur@math.northwestern.edu}}
\asciiemail{dominic.joyce@lincoln.oxford.ac.uk, salur@math.northwestern.edu}

\begin{abstract} McLean proved that the moduli space of
coassociative deformations of a compact coassociative
4--submanifold $C$ in a $G_2$--manifold $(M,\vp,g)$ is a smooth
manifold of dimension equal to $b^2_+(C)$. In this paper, we show
that the moduli space of coassociative deformations of a
noncompact, {\it asymptotically cylindrical\/} coassociative
4--fold $C$ in an asymptotically cylindrical $G_2$--manifold
$(M,\vp,g)$ is also a smooth manifold. Its dimension is the
dimension of the positive subspace of the image of $H^2_{\rm
cs}(C,\R)$ in~$H^2(C,\R)$.
\end{abstract}

\asciiabstract{%
McLean proved that the moduli space of coassociative deformations of a
compact coassociative 4-submanifold C in a G_2-manifold (M,phi,g) is a
smooth manifold of dimension equal to b^2_+(C).  In this paper, we
show that the moduli space of coassociative deformations of a
noncompact, asymptotically cylindrical coassociative 4-fold C in an
asymptotically cylindrical G_2-manifold (M,phi,g) is also a smooth
manifold.  Its dimension is the dimension of the positive subspace of
the image of H^2_cs(C,R) in H^2(C,R).}


\maketitlepage

\section{Introduction}
\label{co1}

Let $(M,g)$ be a Riemannian 7--manifold whose holonomy
group Hol$(g)$ is the exceptional holonomy group $G_2$
(or, more generally, a subgroup of $G_2$). Then $M$ is
naturally equipped with a constant 3--form $\vp$ and
4--form $*\vp$. We call $(M,\vp,g)$ a $G_2$--{\it manifold}.
{\it Complete} examples of Riemannian 7--manifolds with
holonomy $G_2$ were constructed by Bryant and Salamon
\cite{BrSa}, and {\it compact\/} examples by Joyce
\cite{Joyc1} and Kovalev~\cite{Kova}.

Now $\vp$ and $*\vp$ are {\it calibrations} on $M$, in the sense
of Harvey and Lawson \cite{HaLa}. The corresponding calibrated
submanifolds in $M$ are called {\it associative $3$--folds} and
{\it coassociative $4$--folds}, respectively. They are distinguished
classes of minimal 3-- and 4--submanifolds in $(M,g)$ with a rich
structure, that can be thought of as analogous to complex curves
and surfaces in a Calabi--Yau 3--fold.

Harvey and Lawson \cite{HaLa} introduced four types of calibrated
geometries. Special Lagrangian submanifolds of Calabi--Yau
manifolds, associative and coassociative submanifolds of $G_2$
manifolds and Cayley submanifolds of Spin(7) manifolds. Calibrated
geometries have been of growing interest over the past few years
and represent one of the most mysterious classes of minimal
submanifolds \cite{Leun1}, \cite{Leun2}. A great deal of progress
has been made recently in the field of special Lagrangian
submanifolds that arise in mirror symmetry for Calabi--Yau
manifolds and plays a significant role in string theory, for
references see \cite{Joy1}. As one might expect, another promising
direction for future investigation is calibrated submanifolds in
$G_2$ and Spin(7) manifolds. Recently, some progress has been made
in constructing such submanifolds \cite{IKO,Lot1,Lot2} and in
understanding their deformations \cite{AkSa,Leun3}.

The {\it deformation theory} of {\it compact\/} calibrated
submanifolds was studied by McLean \cite{McLe}. He showed
that if $C$ is a compact coassociative $4$--fold in a
$G_2$--manifold $(M,\vp,g)$, then the moduli space $\M_C$ of
coassociative deformations of $C$ is smooth, with
dimension~$b^2_+(C)$.

This paper proves an analogue of McLean's theorem for a special
class of {\it noncompact\/} coassociative 4--folds. The situation
we are interested in is when $(M,\vp,g)$ is an {\it asymptotically
cylindrical\/} $G_2$--manifold, that is, it is a noncompact
7--manifold with one end asymptotic to the cylinder $X\t\R$ on a
{\it Calabi--Yau $3$--fold\/} $X$. The natural class of noncompact
coassociative 4--folds in $M$ are {\it asymptotically
cylindrical\/} coassociative 4--folds $C$, asymptotic at infinity
in $M$ to a cylinder $L\t\R$, where $L$ is a {\it special
Lagrangian $3$--fold\/} in $X$, with phase~$i$. Understanding the
deformations of such submanifolds when the ambient $G_2$--manifold
decomposes into connected sum of two pieces will provide the
necessary technical framework towards completing the Floer
homology program for coassociative submanifolds, \cite{Leun2}.

In particular, we prove the following theorem.

\begin{thm} Let\/ $(M,\vp,g)$ be a $G_2$--manifold asymptotic to
$X\t(R,\iy)$ with decay rate $\al<0$, where $X$ is a Calabi--Yau
$3$--fold. Let\/ $C$ be a coassociative $4$--fold in $M$ asymptotic
to $L\t(R',\iy)$ for $R'>R$ with decay rate $\be$ for $\al\le\be<0$,
where $L$ is a special Lagrangian $3$--fold in $X$ with phase~$i$.

If $\ga<0$ is small enough then the moduli space $\M_C^\ga$ of
asymptotically cylindrical coassociative submanifolds in $M$ close
to $C$, and asymptotic to $L\t(R',\iy)$ with decay rate $\ga$,
is a smooth manifold of dimension $\dim V_+$, where $V_+$ is the
positive subspace of the image of\/ $H^2_{\rm cs}(C,\R)$ in~$H^2(C,\R)$.
\label{co1thm}
\end{thm}

The principal analytic tool we shall use to prove this is
the theory of {\it weighted Sobolev spaces\/} on manifolds
with ends, developed by Lockhart and McOwen \cite{Lock,LoMc}.
The important fact is that elliptic partial differential
operators on exterior forms such as $\d+\d^*$ or $\d^*\d+\d\d^*$
on the noncompact 4--manifold $C$ are {\it Fredholm operators}
between appropriate Banach spaces of forms, and we can
describe their kernels and cokernels.

Results similar to Theorem \ref{co1thm} on the deformations
of classes of noncompact {\it special Lagrangian} $m$--folds
were proved by Marshall \cite{Mars} and Pacini \cite{Paci} for
{\it asymptotically conical\/} special Lagrangian $m$--folds,
and by Joyce \cite{Joyc1,Joyc2} for special Lagrangian $m$--folds
with {\it isolated conical singularities}. Marshall and Joyce
also use the Lockhart--McOwen framework, but Pacini uses a
different analytical approach due to Melrose \cite{Mel1,Mel2}.
Note also that Kovalev \cite{Kova} constructs compact
$G_2$--manifolds by gluing together two noncompact, asymptotically
cylindrical $G_2$--manifolds.

We begin in Section~\ref{co2} with an introduction to $G_2$--manifolds
and coassociative submanifolds, including a sketch of the
proof of McLean's theorem on deformations of compact
coassociative 4--folds, and the definitions of {\it
asymptotically cylindrical\/} $G_2$--manifolds and
coassociative 4--folds. Section \ref{co3} introduces
the weighted Sobolev spaces of Lockhart and McOwen,
and determines the kernel and cokernel of the
elliptic operator $\d_++\d^*$ on $C$ used in the proof.
Finally, Section~\ref{co4} proves Theorem \ref{co1thm}, using
Banach space techniques and elliptic regularity.

\begin{rem} In \cite{Kova}, Kovalev constructs asymptotically
cylindrical manifolds $X$ with holonomy $SU(3)$. Then
$X\times{\mathcal S}^1$ is an asymptotically cylindrical
$G_2$--manifold, though with holonomy $SU(3)$ rather than
$G_2$. One can find examples of asymptotically cylindrical
coassociative 4--folds $C$ in $X\times{\mathcal S}^1$ of two
types:
\begin{itemize}
\item[(a)] $C=C'\times{\rm pt}$, for $C'$ an asymptotically
cylindrical complex surface in $X$; or
\item[(b)] $C=L\times{\mathcal S}^1$, for $L$ an
asymptotically cylindrical special Lagrangian 3--fold in $X$,
with phase $i$.
\end{itemize}

Examples of type (a) can be constructed using algebraic
geometry: if $X=\bar X\sm D$ for $\bar X$ a Fano 3--fold and $D$ a
smooth divisor in $\bar X$, then we can take $C=\bar C\sm D$ for
$\bar C$ a smooth divisor in $\bar X$ intersecting $D$ transversely.
Examples of type (b) can be found by choosing the Calabi--Yau
3--fold $(X,J,\om,\Om)$ to have an {\it antiholomorphic involution}
$\si\co X\ra X$ with $\si^*(J)=-J$, $\si^*(\om)=-\om$ and $\si^*
(\Om)=-\bar\Om$. Then the fixed points $L$ of $\si$ are a
special Lagrangian 3--fold with phase $i$, and each infinite
end of $L$ is asymptotically cylindrical.

We can then apply Theorem \ref{co1thm} to these examples.
One can show that if $\ti C$ is a small deformation of
a coassociative 4--fold $C$ of type (a) or (b) then $\ti C$
is also of type (a) or (b) and thus, Theorem \ref{co1thm}
implies analogous results on the deformation
theory of asymptotically cylindrical complex surfaces and
special Lagrangian 3--folds in asymptotically cylindrical
Calabi--Yau 3--folds.
\end{rem}


\section{Introduction to $G_2$ geometry}
\label{co2}

We now give background material on $G_2$--manifolds and
their coassociative submanifolds that will be needed
later. A good reference on $G_2$ geometry is Joyce \cite[Sections
10--12]{Joyc1}, and a good reference on calibrated geometry is Harvey
and Lawson~\cite{HaLa}.

\subsection{$G_2$--manifolds and coassociative submanifolds}
\label{co21}

Let $(x_1,\dots,x_7)$ be coordinates on $\R^7$. Write
$\d{\bf x}_{ij\ldots l}$ for the exterior form $\d x_i\w\d x_j
\w\cdots\w\d x_l$ on $\R^7$. Define a metric $g_0$, a 3--form $\vp_0$
and a 4--form $*\vp_0$ on $\R^7$ by~$g_0=\d x_1^2+\cdots+\d x_7^2$,
\begin{equation}
\begin{split}
\vp_0&=\d{\bf x}_{123}+\d{\bf x}_{145}+\d{\bf x}_{167}
+\d{\bf x}_{246}-\d{\bf x}_{257}-\d{\bf x}_{347}-\d{\bf x}_{356}
\;\>\text{and}\\
*\vp_0&=\d{\bf x}_{4567}+\d{\bf x}_{2367}+\d{\bf x}_{2345}
+\d{\bf x}_{1357}-\d{\bf x}_{1346}-\d{\bf x}_{1256}-\d{\bf x}_{1247}.
\end{split}
\label{co2eq1}
\end{equation}
The subgroup of $GL(7,\R)$ preserving $\vp_0$ is the {\it exceptional
Lie group} $G_2$. It also preserves $g_0,*\vp_0$ and the orientation
on $\R^7$. It is a compact, semisimple, 14--dimensional
Lie group, a subgroup of~$\SO(7)$.

A {\it $G_2$--structure} on a 7--manifold $M$ is a principal
subbundle of the frame bundle of $M$, with structure group $G_2$.
Each $G_2$--structure gives rise to a 3--form $\vp$ and a metric $g$
on $M$, such that every tangent space of $M$ admits an isomorphism
with $\R^7$ identifying $\vp$ and $g$ with $\vp_0$ and $g_0$
respectively. By an abuse of notation, we will refer to $(\vp,g)$
as a $G_2$--structure.

\begin{prop} Let\/ $M$ be a $7$--manifold and\/ $(\vp,g)$ a
$G_2$--structure on $M$. Then the following are equivalent:
\begin{itemize}
\setlength{\parsep}{0pt}
\setlength{\itemsep}{0pt}
\item[{\rm(i)}] $\Hol(g)\subseteq G_2$, and\/ $\vp$ is the
induced\/ $3$--form,
\item[{\rm(ii)}] $\na\vp=0$ on $M$, where $\na$ is the
Levi--Civita connection of $g$, and
\item[{\rm(iii)}] $\d\vp=\d^*\vp=0$ on $M$.
\end{itemize}
\label{co2prop1}
\end{prop}

We call $\na\vp$ the {\it torsion} of the $G_2$--structure
$(\vp,g)$, and when $\na\vp=0$ the $G_2$--structure is
{\it torsion-free}. A triple $(M,\vp,g)$ is called a
$G_2$--{\it manifold} if $M$ is a 7--manifold and $(\vp,g)$ a
torsion-free $G_2$--structure on $M$. If $g$ has holonomy
$\Hol(g)\subseteq G_2$, then $g$ is Ricci-flat. For explicit,
{\it complete} examples of $G_2$--manifolds see Bryant and
Salamon \cite{BrSa}, and for {\it compact\/} examples
see Joyce \cite{Joyc1} and Kovalev \cite{Kova}. Here are
the basic definitions in {\it calibrated geometry}, due
to Harvey and Lawson~\cite{HaLa}.

\begin{dfn} Let $(M,g)$ be a Riemannian manifold. An {\it oriented
tangent\/ $k$--plane} $V$ on $M$ is a vector subspace $V$ of
some tangent space $T_xM$ to $M$ with $\dim V=k$, equipped
with an orientation. If $V$ is an oriented tangent $k$--plane
on $M$ then $g\vert_V$ is a Euclidean metric on $V$, so
combining $g\vert_V$ with the orientation on $V$ gives a
natural {\it volume form} $\vol_V$ on $V$, which is a
$k$--form on~$V$.

Now let $\vp$ be a closed $k$--form on $M$. We say that
$\vp$ is a {\it calibration} on $M$ if for every oriented
$k$--plane $V$ on $M$ we have $\vp\vert_V\le \vol_V$. Here
$\vp\vert_V=\al\cdot\vol_V$ for some $\al\in\R$, and
$\vp\vert_V\le\vol_V$ if $\al\le 1$. Let $N$ be an
oriented submanifold of $M$ with dimension $k$. Then
each tangent space $T_xN$ for $x\in N$ is an oriented
tangent $k$--plane. We call $N$ a {\it calibrated
submanifold\/} if $\vp\vert_{T_xN}\!=\!\vol_{T_xN}$
for all~$x\!\in\!N$.
\label{co2def1}
\end{dfn}

Calibrated submanifolds are automatically {\it minimal
submanifolds} (see \cite[Theorem~II.4.2]{HaLa}). There are two
natural classes of calibrated submanifolds in $G_2$--manifolds.

\begin{dfn} Let $(M,\vp,g)$ be a $G_2$--{\it manifold}, as
above. Then the 3--form $\vp$ is a {\it calibration} on
$(M,g)$. We define an {\it associative $3$--fold} in $M$ to
be a 3--submanifold of $M$ calibrated with respect to $\vp$.
Similarly, the Hodge star $*\vp$ of $\vp$ is a calibration
4--form on $(M,g)$. We define a {\it coassociative $4$--fold}
in $M$ to be a 4--submanifold of $M$ calibrated with respect
to~$*\vp$.
\label{co2def2}
\end{dfn}

McLean \cite[Prop.~4.4]{McLe} gives an alternative
description of coassociative 4--folds:

\begin{prop} Let\/ $(M,\vp,g)$ be a $G_2$--manifold, and\/ $C$
a $4$--dimensional submanifold of\/ $M$. Then $C$ admits an
orientation making it into a coassociative $4$--fold if and
only if\/~$\vp\vert_C\equiv 0$.
\label{co2prop2}
\end{prop}

\subsection{Deformations of compact coassociative $4$--folds}
\label{co22}

Here is the main result in the {\it deformation theory} of
coassociative 4--folds, proved by McLean \cite[Theorem~4.5]{McLe}.
As our sign conventions for $\vp_0,*\vp_0$ in \eq{co2eq1} are
different to McLean's, we use self-dual 2--forms in place of
McLean's anti-self-dual 2--forms.

\begin{thm} Let\/ $(M,\vp,g)$ be a $G_2$--manifold, and\/ $C$ a
compact coassociative $4$--fold in $M$. Then the moduli space
$\M_C$ of coassociative $4$--folds isotopic to $C$ in
$M$ is a smooth manifold of dimension~$b^2_+(C)$.
\label{co2thm}
\end{thm}

\begin{proof}[Sketch proof] Suppose for simplicity that $C$ is an
embedded submanifold. There is a natural orthogonal decomposition
$TM\vert_C=TC\op\nu$, where $\nu\ra C$ is the {\it normal bundle}
of $C$ in $M$. There is a natural isomorphism
$\nu\cong\La^2_+T^*C$, constructed as follows. Let $x\in C$ and
$V\in\nu_x$. Then $V$ lies in $T_xM$, so
$V\cdot\vp\vert_x\in\La^2T_x^*M$, and
$(V\cdot\vp\vert_x)\vert_{T_xC}\in\La^2T_x^*C$. Moreover
$(V\cdot\vp\vert_x)\vert_{T_xC}$ actually lies in $\La^2_+T_x^*C$,
the bundle of {\it self-dual\/ $2$--forms} on $C$, and the map
$V\mapsto(V\cdot\vp\vert_x)\vert_{T_xC}$ defines an {\it
isomorphism}~$\nu\,\smash{{\buildrel\cong\over\longra}}
\,\La^2_+T^*C$.

For small $\ep>0$, write $B_\ep(\nu)$ for the subbundle of
$\nu$ with fibre at $x$ the open ball about 0 in $\nu\vert_x$
with radius $\ep$. Then the exponential map $\exp\co \nu\ra M$
induces a diffeomorphism between $B_\ep(\nu)$ and a small
{\it tubular neighbourhood\/} $T_C$ of $C$ in $M$. The
isomorphism $\nu\cong\La^2_+T^*C$ gives a diffeomorphism
$\exp\co B_\ep(\La^2_+T^*C)\ra T_C$. Let $\pi\co T_C\ra C$ be
the obvious projection.

Under this identification, submanifolds $\ti C$ in $T_C\subset M$
which are $C^1$ close to $C$ are identified with the {\it graphs}
$\Ga_{\ze^2_+}$ of small smooth sections $\ze^2_+$ of $\La^2_+T^*C$
lying in $B_\ep(\La^2_+T^*C)$. For each $\ze^2_+\in C^\iy\bigl(
B_\ep(\La^2_+T^*C)\bigr)$ the graph $\Ga_{\ze^2_+}$ is a
4--submanifold of $B_\ep(\La^2_+T^*C)$, and so $\ti C=\exp
(\Ga_{\ze^2_+})$ is a 4--submanifold of $T_C$. We need to know:
which 2--forms $\ze^2_+$ correspond to {\it coassociative} 4--folds
$\ti C$ in~$T_C$?

$\ti C$ is coassociative if $\vp\vert_{\ti C}\equiv 0$. Now
$\pi\vert_{\ti C}\co \ti C\ra C$ is a diffeomorphism, so we can push
$\vp\vert_{\ti C}$ down to $C$, and regard it as a function of
$\ze^2_+$. That is, we define
\begin{equation*}
Q\co C^\iy\bigl(B_\ep(\La^2_+T^*C)\bigr)\ra C^\iy(\La^3T^*C)
\;\>\text{by}\;\>
Q(\ze^2_+)=\pi_*(\vp\vert_{\exp(\Ga_{\ze^2_+})}).
\end{equation*}
Then the moduli space $\M_C$ is locally isomorphic near $C$ to the
set of small self-dual 2--forms $\ze^2_+$ on $C$ with
$\vp\vert_{\exp(\Ga_{\ze^2_+})}\equiv 0$, that is, to a
neighborhood of 0 in~$Q^{-1}(0)$.

To understand the equation $Q(\ze^2_+)=0$, note that at $x\in C$,
$Q(\ze^2_+)\vert_x$ depends on the tangent space to $\Ga_{\ze^2_+}$
at $\ze^2_+\vert_x$, and so on $\ze^2_+\vert_x$ and
$\na\ze^2_+\vert_x$. Thus the functional form of $Q$ is
\begin{equation*}
Q(\ze^2_+)\vert_x=F\bigl(x,\ze^2_+\vert_x,\na\ze^2_+\vert_x\bigr)
\quad\text{for $x\in C$,}
\end{equation*}
where $F$ is a smooth function of its arguments. Hence
$Q(\ze^2_+)=0$ is a {\it nonlinear first order PDE} in
$\ze^2_+$. As $\vp$ is closed, $\vp\vert_C\equiv 0$, and
$\Ga_{\ze^2_+}$ is isotopic to $C$, we see that
$\vp\vert_{\Ga_{\ze^2_+}}$ is an exact 3--form on $\Ga_{\ze^2_+}$,
so that $Q(\ze^2_+)$ is {\it exact}. The {\it linearization} $\d
Q(0)$ of $Q$ at $\ze^2_+=0$ is
\begin{equation*}
\d Q(0)(\be)=\lim_{\ep\ra 0}\bigl(\ep^{-1}Q(\ep\be)\bigr)=\d\be.
\end{equation*}
Therefore $\Ker(\d Q(0))$ is the vector space $\H^2_+$
of {\it closed self-dual\/ $2$--forms} $\be$ on $C$, which by Hodge
theory is a finite-dimensional vector space isomorphic to
$H^2_+(C,\R)$, with dimension $b^2_+(C)$. This is the {\it Zariski
tangent space} of $\M_C$ at $C$, the {\it infinitesimal
deformation space} of $C$ as a coassociative 4--fold.

To complete the proof we must show that $\M_C$ is
locally isomorphic to its Zariski tangent space ${\mathcal
H}^2_+$, and so is a smooth manifold of dimension $b^2_+(C)$. To
do this rigorously requires some technical analytic machinery,
which is passed over in a few lines in \cite[p.~731]{McLe}. Here
is one way to do it.

As $Q$ maps from $\La^2_+T^*C$ with fibre $\R^3$ to $\La^3T^*C$
with fibre $\R^4$, it is {\it overdetermined}, and not {\it
elliptic}. To turn it into an elliptic operator, define
\begin{equation}
\begin{gathered}
P\co C^\iy\bigl(B_\ep(\La^2_+T^*C)\bigr)\t C^\iy(\La^4T^*C)\ra
C^\iy(\La^3T^*C) \\
\text{by}\qquad P(\ze^2_+,\ze^4)=Q(\ze^2_+)+\d^*\ze^4.
\end{gathered}
\label{co2eq2}
\end{equation}
Then the linearization of $P$ at $(0,0)$ is
\begin{equation*}
\d P(0,0)\co (\ze^2_+,\ze^4)\mapsto \d\ze^2_++\d^*\ze^4,
\end{equation*}
which is elliptic. Since ellipticity is an open condition, $P$ is
elliptic near $(0,0)$ in~$C^\iy\bigl(B_\ep(\La^2_+T^*C)\bigr)\t
C^\iy(\La^4T^*C)$.

Suppose $P(\ze^2_+,\ze^4)=0$. Then $Q(\ze^2_+)=-\d^*\ze^4$, so
\begin{equation*}
\lnm{\d^*\ze^4}{2}^2=-\bigl\langle\d^*\ze^4,Q(\ze^2_+)\bigr\rangle_{L^2}=
-\bigl\langle\ze^4,\d(Q(\ze^2_+))\bigr\rangle_{L^2}=0,
\end{equation*}
by integration by parts, since $Q(\ze^2_+)$ is {\it exact}. Hence
$P(\ze^2_+,\ze^4)=0$ if and only if $Q(\ze^2_+)=\d^*\ze^4=0$. But
4--forms with $\d^*\ze^4=0$ are constant, and the vector space of
such $\ze^4$ is $H^4(C,\R)$. Thus,~$P^{-1}(0)=Q^{-1}(0)\t H^4(C,\R)$.

Because $C^\iy(\La^2_+T^*C),C^\iy(\La^3T^*C)$ are not {\it Banach
spaces}, we extend $P$ in \eq{co2eq2} to act on {\it Sobolev
spaces} $L^p_{l+2}(\La^2_+T^*C),L^p_{l+2}(\La^3T^*C)$ for $p>4$
and $l\ge 0$, giving
\begin{equation*}
\begin{split}
&\ti P\co L^p_{l+2}\bigl(B_\ep(\La^2_+T^*C)\bigr)\t
L^p_{l+2}(\La^4T^*C)\ra L^p_{l+1}(\La^3T^*C),\\
&\ti P\co (\ze^2_+,\ze^4)\mapsto
\pi_*(\vp\vert_{\Ga_{\ze^2_+}})+\d^*\ze^4.
\end{split}
\end{equation*}
Then $\ti P$ is a {\it smooth map of Banach manifolds}.

Let $\H^3$ be the vector space of closed and coclosed
3--forms on $C$, so that $\H^3\cong H^3(C,\R)$ by Hodge
theory, and $V^p_{l+1}$ be the Banach subspace of
$L^p_{l+1}(\La^3T^*C)$ $L^2$--orthogonal to $\H^3$. Then one can
show that $\ti P$ maps into $V^p_{l+1}$, and the linearization
\begin{align*}
\d\ti P(0,0)&\co L^p_{l+2}(\La^2_+T^*C)\t L^p_{l+2}(\La^4T^*C)\ra
V^p_{l+1},\\
\d\ti P(0,0)&\co (\ze^2_+,\ze^4)\mapsto \d\ze^2_++\d^*\ze^4
\end{align*}
is then {\it surjective} as a map of Banach spaces.

Thus, $\ti P\co L^p_{l+2}\bigl(B_\ep(\La^2_+T^*C)\bigr)\t
L^p_{l+2}(\La^4T^*C)\ra V^p_{l+1}$ is a smooth map of Banach
manifolds, with $\d\ti P(0,0)$ surjective. The {\it Implicit
Mapping Theorem for Banach spaces} (Theorem \ref{co4thm}) now
implies that $\ti P^{-1}(0)$ is near 0 a smooth submanifold,
locally isomorphic to $\Ker(\d\ti P(0))$. But $\ti P(\ze^2_+,
\ze^4)=0$ is an {\it elliptic equation} for small $\ze^2_+,\ze^4$,
and so elliptic regularity implies that solutions $(\ze^2_+,\ze^4)$
are smooth.

Therefore $\ti P^{-1}(0)=P^{-1}(0)$ near 0, and also $\Ker(\d\ti
P(0,0))=\Ker(\d P(0,0))$. Hence $P^{-1}(0)$ is, near $(0,0)$, a
smooth manifold locally isomorphic to the kernel $\Ker(\d P(0,0))$. So from
above $Q^{-1}(0)$ is near 0 a smooth manifold locally isomorphic
to $\Ker(\d Q(0))$. Thus, $\M_C$ is near $C$ a smooth manifold
locally isomorphic to $H^2_+(C,\R)$. This completes the proof.
\end{proof}

\subsection{Asymptotically cylindrical $G_2$--manifolds and
coassociative 4--folds}
\label{co23}

We first define {\it cylindrical\/} and {\it asymptotically
cylindrical\/} $G_2$--manifolds.

\begin{dfn} A $G_2$--manifold $(M_0,\vp_0,g_0)$ is called
{\it cylindrical\/} if $M_0=X\t\R$ and $(\vp_0,g_0)$ is
compatible with this product structure, that is,
\begin{equation*}
\vp_0=\Re\Om+\om\w\d t \qquad\text{and}\qquad g_0=g_X+\d t^2,
\end{equation*}
where $X$ is a (connected, compact) Calabi--Yau 3--fold with
K\"{a}hler form $\om$, Riemannian metric $g_X$ and holomorphic
(3,0)-form~$\Om$.
\label{co2def3}
\end{dfn}

\begin{dfn} A connected, complete $G_2$--manifold $(M,\vp,g)$ is
called {\it asymptotically cylindrical with decay rate\/} $\al$
if there exists a cylindrical $G_2$--manifold $(M_0,\vp_0,g_0)$
with $M_0=X\t\R$ as above, a compact subset $K\subset M$,
a real number $R$, and a diffeomorphism $\Psi\co X\t(R,\iy)\ra
M\sm K$ such that $\Psi^*(\vp)=\vp_0+\d\xi$ for some smooth
2--form $\xi$ on $X\t(R,\iy)$ with $\bmd{\na^k\xi}=O(e^{\al t})$
on $X\t(R,\iy)$ for all $k\ge 0$, where $\na$ is the Levi--Civita
connection of the cylindrical metric~$g_0$.
\label{co2def4}
\end{dfn}

The point of this is that $M$ has one end modelled
on $X\t(R,\iy)$, and as $t\ra\iy$ in $(R,\iy)$ the
$G_2$--structure $(\vp,g)$ on $M$ converges to order
$O(e^{\al t})$ to the cylindrical $G_2$--structure
on $X\t(R,\iy)$, with all of its derivatives. We
suppose $M$ and $X$ are connected, that is, we allow
$M$ to have {\it only one end}. 

This is because one can use Cheeger--Gromoll splitting theorem \cite{CG}
to show that an orientable, connected, asymptotically cylindrical
Riemannian manifold with Ricci-flat metric $g$ can have at most 2
cylindrical ends. In the case when there are 2 cylindrical ends then 
there is reduction in the holonomy group Hol$(g)$ and $(M,g)$ is a cylinder.
One can also show that reduction in holonomy can be obtained by just using
the analytic set-up for Fredholm properties of an elliptic
operator on noncompact manifolds, \cite{Salu01}. 



Here are the analogous definitions for coassociative submanifolds.

\begin{dfn} Let $(M_0,\vp_0,g_0)$ and $X$ be as in Definition
\ref{co2def3}. A submanifold $C_0$ of $M_0$ is called {\it
cylindrical\/} if $C_0=L\t\R$ for some compact submanifold $L$ in
$X$, not necessarily connected. $C_0$ is {\it coassociative} if
and only if $L$ is a {\it special Lagrangian $3$--fold with phase}
$i$ in the Calabi--Yau 3--fold~$X$. \label{co2def5}
\end{dfn}

\begin{dfn} Let $(M_0,\vp_0,g_0)$, $X$, $(M,\vp,g)$, $K,\Psi$
and $\al$ be as in Definitions \ref{co2def3} and \ref{co2def4},
and let $C_0=L\t\R$ be a cylindrical coassociative 4--fold
in $M_0$, as in Definition~\ref{co2def5}.

A connected, complete coassociative 4--fold $C$ in $(M,\vp,g)$ is
called {\it asymptotically cylindrical with decay rate $\be$} for
$\al\le\be<0$ if there exists a compact subset $K'\subset C$,
a normal vector field $v$ on $L\t(R',\iy)$ for some $R'>R$,
and a diffeomorphism $\Phi\co L\t(R',\iy)\ra C\sm K'$ such that
the diagram
\begin{equation}
\begin{gathered}
\xymatrix{X\t (R',\iy) \ar[d]^\subset & L\t (R',\iy)
\ar[l]^{\exp_v} \ar[r]_{\Phi} & (C\sm K')
\ar[d]^\subset \\
X\t (R,\iy)\ar[rr]^\Psi && (M\sm K) }
\end{gathered}
\label{co2eq3}
\end{equation}
commutes, and $\bmd{\na^kv}=O(e^{\be t})$ on $L\t(R',\iy)$ for all~$k\ge 0$.
\label{co2def6}
\end{dfn}

Here we require $C$ but {\it not\/} $L$ to be connected, that is,
we allow $C$ to have {\it multiple ends}. The point of Definition
\ref{co2def6} is to find a good way to say that a submanifold
$C$ in $M$ is asymptotic to the cylinder $C_0$ in $M_0=X\t\R$
as $t\ra\iy$ in $\R$, to order $O(e^{\be t})$. We do this by
writing $C$ near infinity as the graph of a {\it normal vector
field\/} $v$ to $C_0=L\t\R$ in $M_0=X\t\R$, and requiring
$v$ and its derivatives to be~$O(e^{\be t})$.

\section{Infinitesimal deformations of $C$}
\label{co3}

Let $(M,\vp,g)$ be an asymptotically cylindrical $G_2$--manifold
asymptotic to $X\t(R,\iy)$, and $C$ an asymptotically cylindrical
coassociative 4--fold in $M$ asymptotic to $L\t(R',\iy)$. We
wish to study the moduli space $\M_C^\ga$ of asymptotically
cylindrical deformations $\ti C$ of $C$ in $M$ with rate $\ga$.
To do this we modify the proof of Theorem \ref{co2thm} in Section~\ref{co2},
for the case when $C$ is compact. There we modelled $\M_C$ on
$\ti P^{-1}(0)$ for a nonlinear map $\ti P$ between Banach spaces,
whose linearization $\d\ti P(0,0)$ at 0 was the Fredholm map
between Sobolev spaces
\begin{equation}
\d_++\d^*\co L^p_{l+2}(\La^2_+T^*C)\t L^p_{l+2}(\La^4T^*C)\longra
L^p_{l+1}(\La^3T^*C).
\label{co3eq1}
\end{equation}
Now when $C$ is not compact, as in the asymptotically cylindrical
case, \eq{co3eq1} is not in general Fredholm, and the proof of
Theorem \ref{co2thm} fails. To repair it we use the analytical
framework for asymptotically cylindrical manifolds developed
by Lockhart and McOwen in \cite{Lock,LoMc}, involving
{\it weighted Sobolev spaces} $L^p_{k,\ga}(\La^rT^*C)$. Roughly
speaking, elements of $L^p_{k,\ga}(\La^rT^*C)$ are $L^p_k$
$r$--forms on $C$ which decay like $O(e^{\ga t})$ on the
end $L\t(R',\iy)$. This has the advantage of building the
decay rate $\ga$ into the proof from the outset.

This section will study the weighted analogue of \eq{co3eq1},
\begin{equation}
\d_++\d^*\co L^p_{l+2,\ga}(\La^2_+T^*C)\t L^p_{l+2,\ga}(\La^4T^*C)
\longra L^p_{l+1,\ga}(\La^3T^*C),
\label{co3eq2}
\end{equation}
for small $\ga<0$. It will be shown in Section~\ref{co4} to be the
linearization at 0 of a nonlinear operator $P$ for which
$\M_C^\ga$ is locally modelled on~$P^{-1}(0)$.

Section \ref{co31} introduces weighted Sobolev spaces, and
the Lockhart--McOwen theory of elliptic operators between
them. Then Sections \ref{co32} and \ref{co33} compute the
{\it kernel\/} and {\it cokernel\/} of \eq{co3eq2} for
small $\ga<0$, and Section~\ref{co34} determines the set of
rates $\ga$ for which \eq{co3eq2} is {\it Fredholm}.

\subsection{Elliptic operators on asymptotically cylindrical manifolds}
\label{co31}

We now sketch parts of the theory of analysis on manifolds
with cylindrical ends due to Lockhart and McOwen \cite{Lock,LoMc}. We
begin with some elementary definitions.

\begin{dfn} Let $(C,g)$ be an {\it asymptotically cylindrical
Riemannian manifold}. That is, there is a Riemannian
cylinder $(L\t\R,g_0)$ with $L$ compact, a compact subset
$K'\subset C$ and a diffeomorphism $\Phi\co C\sm K'\ra L\t(R',\iy)$
such that
\begin{equation*}
\na_0^k\bigl(\Phi_*(g)-g_0\bigr)=O(e^{\be t})
\quad\text{for all $k\ge 0$}
\end{equation*}
for some rate $\be<0$, where $\na_0$ is the
Levi--Civita connection of $g_0$ on~$L\t\R$.

Let $E_0$ be a {\it cylindrical vector bundle} on
$L\t\R$, that is, a vector bundle on $L\t\R$
invariant under translations in $\R$. Let $h_0$ be
a {\it smooth family of metrics} on the fibres of
$E_0$ and $\na_{\sst E_0}$ a {\it connection} on $E_0$
preserving $h_0$, with $h_0,\na_{\sst E_0}$ invariant
under translations in~$\R$.

Let $E$ be a vector bundle on $C$ equipped with metrics $h$ on the
fibres, and a connection $\na_{\sst E}$ on $E$ preserving $h$. We
say that $E,h,\na_{\sst E}$ are {\it asymptotic to}
$E_0,h_0,\na_{\sst E_0}$ if there exists an identification
$\Phi_*(E)\cong E_0$ on $L\t(R',\iy)$ such that
$\Phi_*(h)=h_0+O(e^{\be t})$ and $\Phi_*( \na_{\sst E})=\na_{\sst
E_0}+O(e^{\be t})$ as $t\ra\iy$. Then we call $E,h,\na_{\sst E}$
{\it asymptotically cylindrical}.

Choose a smooth function $\rho\co C\ra\R$ such that $\Phi^*(\rho)
\equiv t$ on $L\t(R',\iy)$. This prescribes $\rho$ on $C\sm K'$,
so we only have to extend $\rho$ over the compact set $K'$. For
$p\geq 1$, $k\geq 0$ and $\ga\in\R$ we define the {\it weighted
Sobolev space} $L^p_{k,\ga}(E)$ to be the set of sections $s$ of
$E$ that are locally integrable and $k$ times weakly differentiable
and for which the norm
\begin{equation}
\nm{s}_{L^p_{k,\ga}}=\Bigl(\sum_{j=0}^{k}\int_C
e^{-\ga\rho}\bmd{\na_{\sst E}^js}^p\d V\Bigr)^{1/p}
\label{co3eq3}
\end{equation}
is finite. Then $L^p_{k,\ga}(E)$ is a Banach space. Since
$\rho$ is uniquely determined except on the compact set $K'$,
different choices of $\rho$ give the same space $L^p_{k,\ga}(E)$,
with equivalent norms.
\label{co3def1}
\end{dfn}

For instance, the $r$--forms $E=\La^rT^*C$ on $C$ with
metric $g$ and the Levi--Civita connection are automatically
asymptotically cylindrical, and if $C$ is an oriented
4--manifold then the self-dual 2--forms $\La^2_+T^*C$ are
also asymptotically cylindrical. We consider {\it partial
differential operators} on asymptotically cylindrical manifolds.

\begin{dfn} In the situation of Definition \ref{co3def1},
suppose $E,F$ are two asymptotically cylindrical vector
bundles on $C$, asymptotic to cylindrical vector bundles
$E_0,F_0$ on $L\t\R$. Let $A_0\co C^\iy(E_0)\ra C^\iy(F_0)$
be a linear partial differential operator of order $k$
which is {\it cylindrical}, that is, invariant under
translations in~$\R$.

Suppose $A\co C^\iy(E)\ra C^\iy(F)$ is a linear partial differential
operator of order $k$ on $C$. We say that $A$ is {\it asymptotic
to} $A_0$ if under the identifications $\Phi_*(E)\cong E_0$,
$\Phi_*(F)\cong F_0$ on $L\t(R',\iy)$ we have
$\Phi_*(A)=A_0+O(e^{\be t})$ as $t\ra\iy$ for $\be<0$. Then we
call $A$ an {\it asymptotically cylindrical\/} operator. It is
easy to show that $A$ extends to bounded linear operators
\begin{equation}
A^p_{k+l,\ga}\co L^p_{k+l,\ga}(E)\longra L^p_{l,\ga}(F)
\label{co3eq4}
\end{equation}
for all $p>1$, $l\ge 0$ and~$\ga\in\R$.
\label{co3def2}
\end{dfn}

Now suppose $A$ is an {\it elliptic} operator. \eq{co3eq4} is {\it
Fredholm} if and only if $\ga$ does not lie in a discrete set
$\D_{A_0}\subset\R$, which we now define.

\begin{dfn} In Definition \ref{co3def2}, suppose $A$ and $A_0$
are {\it elliptic} operators on $C$ and $L\t\R$, so that $E,F$
have the same fibre dimensions. Extend $A_0$ to the
complexifications $A_0\co C^\iy(E_0\ot_\R\C)\ra C^\iy(F_0\ot_\R\C)$.
Define $\D_{A_0}$ to be the set of $\ga\in\R$ such that for some
$\de\in\R$ there exists a nonzero section $s\in C^\iy(E_0\ot_\R\C)$
invariant under translations in $\R$ such that~$A_0(e^{(\ga+i\de)t}s)=0$.
\label{co3def3}
\end{dfn}

Then Lockhart and McOwen prove~\cite[Theorem~1.1]{LoMc}:

\begin{thm} Let\/ $(C,g)$ be a Riemannian manifold asymptotic to $(L\t\R,g_0)$, and\/
$A\co C^\iy(E)\ra C^\iy(F)$ an elliptic
partial differential operator on $C$ of order $k$ between
vector bundles $E,F$ on $C$,
asymptotic to the cylindrical elliptic operator $A_0:
C^\iy(E_0)\ra C^\iy(F_0)$ on $L\t\R$. Define $\D_{A_0}$ as above.

Then $\D_{A_0}$ is a discrete subset of\/ $\R$, and for $p>1$,
$l\ge 0$ and $\ga\in\R$, the extension $A^p_{k+l,\ga}\co L^p_{k+l,\ga}
(E)\ra L^p_{l,\ga}(F)$ is Fredholm if and only if\/~$\ga\notin\D_{A_0}$.
\label{co3thm1}
\end{thm}

Suppose $\ga\notin\D_{A_0}$. Then $A^p_{k+l,\ga}$ is Fredholm,
so its kernel $\Ker(A^p_{k+l,\ga})$ is finite-dimensional. Let
$e\in\Ker(A^p_{k+l,\ga})$. Then by an elliptic regularity result
\cite[Theorem~3.7.2]{Lock} we have $e\in L^p_{k+m,\ga}(E)$ for
all $m\ge 0$. The {\it weighted Sobolev Embedding Theorem}
\cite[Theorem~3.10]{Lock} then implies that $e\in L^r_{k+m,\de}
(E)$ for all $r>1$, $m\ge 0$ and $\de>\ga$, and $e$ is smooth.
But $\Ker(A^p_{k+1,\ga})$ is invariant under small changes of
$\ga$ in $\R\sm\D_{A_0}$, so $e\in L^r_{k+m,\ga}(E)$ for
all $r>1$ and $m\ge 0$. This proves:

\begin{prop} For $\ga\notin\D_{A_0}$ the kernel\/ $\Ker(A^p_{k+l,\ga})$
is independent of\/ $p,l$, and is a finite-dimensional vector space of
smooth sections of\/~$E$.
\label{co3prop1}
\end{prop}

When $\ga\notin\D_{A_0}$, as $A^p_{k+l,\ga}$ is Fredholm the
{\it cokernel\/}
\begin{equation*}
\Coker(A^p_{k+l,\ga})=L^p_{l,\ga}(F)\big/A^p_{k+l,\ga}
\bigl(L^p_{k+l,\ga}(E)\bigr)
\end{equation*}
of $A^p_{k+l,\ga}$ is also finite-dimensional. To understand it,
consider the {\it formal adjoint\/} $A^*\co C^\iy(F)\ra C^\iy(E)$ of $A$.
This is also an asymptotically cylindrical linear elliptic partial
differential operator of order $k$ on $C$, with the property that
\begin{equation*}
\langle Ae,f\rangle_{L^2(F)}=\langle e,A^*f\rangle_{L^2(E)}
\end{equation*}
for compactly-supported $e\in C^\iy(E)$ and  $f\in C^\iy(F)$.

Then for $p>1$, $l\ge 0$ and $\ga\notin\D_{A_0}$, the {\it
dual operator} of \eq{co3eq4} is
\begin{equation}
(A^*)^q_{-l,-\ga}\co L^q_{-l,-\ga}(F)\longra L^q_{-k-l,-\ga}(E),
\label{co3eq5}
\end{equation}
where $q>1$ is defined by $\frac{1}{p}+\frac{1}{q}=1$.
Here we mean that $L^q_{-k-l,-\ga}(E)$, $L^q_{-l,-\ga}(F)$
are isomorphic to the Banach space duals of $L^p_{k+l,\ga}(E)$,
$L^p_{l,\ga}(F)$, and these isomorphisms identify
$(A^*)^q_{-l,-\ga}$ with the dual linear map to~\eq{co3eq4}.

Now there is a {\it problem} with \eq{co3eq5}, as it involves
Sobolev spaces with {\it negative numbers of derivatives} $-l,-k-l$.
Such Sobolev spaces exist as spaces of {\it distributions}. But we
can avoid defining or using these spaces, by the following trick.
We are interested in $\Ker\bigl((A^*)^q_{-l,-\ga}\bigr)$, as it is
dual to $\Coker(A^p_{k+l,\ga})$. The elliptic regularity argument
above showing $\Ker(A^p_{k+l,\ga})$ is independent of $l$ also holds
for negative differentiability, so we have $\Ker\bigl((A^*)^q_{-l,-\ga}
\bigr)=\Ker\bigl((A^*)^q_{k+m,-\ga}\bigr)$ for $m\in\Z$, and in
particular for $m\ge 0$. So we deduce:

\begin{prop} In Theorem \ref{co3thm1}, let\/ $A^*$ be the formal
adjoint of\/ $A$. Then for all\/ $\ga\notin\D_{A_0}$, $p,q>1$
with $\frac{1}{p}+\frac{1}{q}=1$ and\/ $l,m\ge 0$ there is a
natural isomorphism
\begin{equation}
\Coker(A^p_{k+l,\ga})\cong\Ker\bigl((A^*)^q_{k+m,-\ga}\bigr)^*.
\label{co3eq6}
\end{equation}
\label{co3prop2}
\end{prop}

When $\ga\notin\D_{A_0}$ we see from \eq{co3eq6} that the
{\it index} of $A^p_{k+l,\ga}$ is
\begin{equation}
\ind(A^p_{k+l,\ga})=\dim\Ker(A^p_{k+l,\ga})
-\dim\Ker\bigl((A^*)^q_{k+m,-\ga}\bigr).
\label{co3eq7}
\end{equation}
Lockhart and McOwen show \cite[Theorem~6.2]{LoMc} that for
$\ga,\de\in\R\sm\D_{A_0}$ with $\ga\le\de$ we have
\begin{equation}
\ind(A^p_{k+l,\de})-\ind(A^p_{k+l,\ga})=
\sum_{\ep\in\D_{A_0}:\ga<\ep<\de}d(\ep),
\label{co3eq8}
\end{equation}
where $d(\ep)\ge 1$ is the dimension of the a vector space
of solutions $s\in C^\iy(E_0\ot_\R\C)$ of a prescribed
form with~$A_0(s)=0$.

\subsection{$\d+\d^*$ and $\d^*\d+\d\d^*$ on an asymptotically
cylindrical manifold}
\label{co32}

Let $(C,g)$ be an oriented asymptotically cylindrical
Riemannian $n$--manifold asymptotic to a Riemannian
cylinder $(L\t\R,g_0)$, where $g_0=g_L+\d t^2$ and
$(L,g_L)$ is a compact oriented Riemannian $(n-1)$--manifold.
Consider the {\it asymptotically cylindrical linear elliptic
operators}
\begin{equation}
\d+\d^*\quad\text{and}\quad \d^*\d+\d\d^*\co \ts\bigop_{k=0}^n
C^\iy(\La^kT^*C)\longra\ts\bigop_{k=0}^nC^\iy(\La^kT^*C).
\label{co3eq9}
\end{equation}
We shall apply the theory of Section~\ref{co31} to study the extensions
\begin{align}
(\d+\d^*)^p_{l+2,\ga}\co
&\ts\bigop_{k=0}^nL^p_{l+2,\ga}(\La^kT^*C)\longra
\ts\bigop_{k=0}^nL^p_{l+1,\ga}(\La^kT^*C),
\label{co3eq10}\\
(\d^*\d+\d\d^*)^p_{l+2,\ga}\co
&\ts\bigop_{k=0}^nL^p_{l+2,\ga}(\La^kT^*C)\longra
\ts\bigop_{k=0}^nL^p_{l,\ga}(\La^kT^*C),
\label{co3eq11}
\end{align}
for $p>1$, $l\ge 0$ and $\ga\in\R$, and their {\it kernels}
and {\it cokernels}.

\begin{lem} We have $\Ker\bigl((\d+\d^*)^p_{l+2,\ga}\bigr)
\subseteq\Ker\bigl((\d^*\d+\d\d^*)^p_{l+2,\ga}\bigr)$ for
all\/ $p>1$, $l\ge 0$ and\/ $\ga\in\R$, and equality holds
if\/~$\ga<0$.
\label{co3lem}
\end{lem}

\begin{proof} Since $\d^*\d+\d\d^*=(\d+\d^*)^2$ we have
$\Ker(\d+\d^*)\subseteq\Ker(\d^*\d+\d\d^*)$ on any space
of twice differentiable forms, giving the inclusion.
Suppose $\ga<0$ and $\chi\in\Ker\bigl((\d^*\d+\d\d^*
)^p_{l+2,\ga}\bigr)$. Write $\chi=\sum_{k=0}^n\chi_k$
for $\chi_k\in L^p_{l+2,\ga}(\La^kT^*C)$. Then $\chi_k
\in\Ker\bigl((\d^*\d+\d\d^*)^p_{l+2,\ga}\bigr)$, as
$\d^*\d+\d\d^*$ takes $k$--forms to $k$--forms.

If $\ga<0$ then each $\chi_k$ lies in $L^2_2(\La^kT^*C)$, and
\begin{equation*}
\lnm{\d\chi_k}{2}^2\!+\!\lnm{\d^*\chi_k}{2}^2\!=\!
\big\langle\d\chi_k,\d\chi_k\big\rangle_{L^2}\!+\!
\big\langle\d^*\chi_k,\d^*\chi_k\big\rangle_{L^2}\!=\!
\big\langle\chi_k,(\d^*\d\!+\!\d\d^*)\chi_k\big\rangle_{L^2}\!=\!0.
\end{equation*}
Thus $\d^*\chi_k=\d\chi_k=0$, so that $\chi_k$ and hence $\chi$
lies in~$\Ker\bigl((\d+\d^*)^p_{l+2,\ga}\bigr)$.
\end{proof}

For $\md{\ga}$ close to zero we can say more about the
kernels of \eq{co3eq10} and~\eq{co3eq11}.

\begin{prop} Suppose $p,q>1$, $l,m\ge 0$ and\/ $\ga<0$ with\/
$\frac{1}{p}+\frac{1}{q}=1$ and\/ $[\ga,-\ga]\cap\D_{(\d+\d^*)_0}
=[\ga,-\ga]\cap\D_{(\d^*\d+\d\d^*)_0}=\{0\}$. Then
\begin{align}
\Ker\bigl((\d+\d^*)^p_{l+2,\ga}\bigr)&=
\Ker\bigl((\d^*\d+\d\d^*)^p_{l+2,\ga}\bigr),
\label{co3eq12}\\
\Ker\bigl((\d+\d^*)^q_{m+2,-\ga}\bigr)&=
\Ker\bigl((\d^*\d+\d\d^*)^q_{m+2,-\ga}\bigr),\qquad\text{and}
\label{co3eq13}\\
\dim\Ker\bigl((\d+\d^*)^q_{m+2,-\ga}\bigr)&=
\dim\Ker\bigl((\d+\d^*)^p_{l+2,\ga}\bigr)+\ts\sum_{k=0}^{n-1}b^k(L).
\label{co3eq14}
\end{align}
Moreover all four kernels consist of smooth closed and coclosed forms.
\label{co3prop3}
\end{prop}

\begin{proof} As $[\ga,-\ga]\cap\D_{(\d+\d^*)_0}=\{0\}$, 
\begin{equation}
\ind\bigl((\d+\d^*)^q_{m+2,-\ga}\bigr)-
\ind\bigl((\d+\d^*)^p_{l+2,\ga}\bigr)
=\ts 2\sum_{k=0}^{n-1}b^k(L).
\label{co3eq15}
\end{equation}
This is because from \eq{co3eq8}, the l.h.s.\ of \eq{co3eq15}
is the dimension of the solution space of $(\d+\d^*)_0\chi=0$ on
$L\t\R$ for $\chi$ independent of $t\in\R$. The space of such
$\chi$ is the direct sum over $k=0,\ldots,n-1$ of the spaces
of $k$--forms $\eta$ and $(k\!+\!1)$--forms $\eta\w\d t$ for
$\eta\in C^\iy(\La^kT^*L)$ with $\d\eta=\d^*\eta=0$. By Hodge
theory we deduce~\eq{co3eq15}.

Now $\d+\d^*$ is {\it formally self adjoint}, that is,
$A^*=A$ in the notation of Section~\ref{co31}. Thus
\begin{gather*}
\ind\bigl((\d+\d^*)^q_{m+2,-\ga}\bigr)=
-\ind\bigl((\d+\d^*)^p_{l+2,\ga}\bigr)=\\
\dim\Ker\bigl((\d+\d^*)^q_{m+2,-\ga}\bigr)
-\dim\Ker\bigl((\d+\d^*)^p_{l+2,\ga}\bigr)
\end{gather*}
by \eq{co3eq7}, and equation \eq{co3eq14} follows
from \eq{co3eq15}. As $[\ga,-\ga]\cap\D_{
(\d^*\d+\d\d^*)_0}=\{0\}$, the same proof shows that
\begin{equation}
\dim\Ker\bigl((\d^*\d+\d\d^*)^q_{m+2,-\ga}\bigr)=
\dim\Ker\bigl((\d^*\d+\d\d^*)^p_{l+2,\ga}\bigr)+
\ts\sum_{k=0}^{n-1}b^k(L),
\label{co3eq16}
\end{equation}
since the solutions of $(\d^*\d+\d\d^*)_0\chi=0$ and
$(\d+\d^*)_0\chi=0$ for $\chi$ on $L\t\R$ independent of
$t$ coincide. Lemma \ref{co3lem} proves \eq{co3eq12}, and
combining this with \eq{co3eq14} and \eq{co3eq16} yields
\begin{equation*}
\dim\Ker\bigl((\d+\d^*)^q_{m+2,-\ga}\bigr)=
\dim\Ker\bigl((\d^*\d+\d\d^*)^q_{m+2,-\ga}\bigr).
\end{equation*}
As the right hand side of \eq{co3eq13} contains the
left by Lemma \ref{co3lem}, this implies~\eq{co3eq13}.

It remains to show the four kernels consist of smooth
closed and coclosed forms. Let $\chi$ lie in one of the
kernels, and write $\chi=\sum_{k=0}^n\chi_k$ for $\chi_k$
a $k$--form. Since $(\d^*\d+\d\d^*)\chi=0$ we have
$(\d^*\d+\d\d^*)\chi_k=0$, as $\d^*\d+\d\d^*$ takes
$k$--forms to $k$--forms. Thus $\chi_k$ lies in the same
kernel, so $(\d+\d^*)\chi_k=0$ by \eq{co3eq12} or
\eq{co3eq13}. But $\d\chi_k$ and $\d^*\chi_k$ lie in
different vector spaces, so $\d\chi_k=\d^*\chi_k=0$ for
all $k$. Hence $\d\chi=\d^*\chi=0$, and $\chi$ is closed
and coclosed. Smoothness follows by elliptic regularity.
\end{proof}

As the forms $\chi$ in $\Ker\bigl((\d+\d^*)^p_{l+2,\ga}\bigr)$
are closed we can map them to de Rham cohomology $H^*(C,\R)$
by $\chi\mapsto[\chi]$. We identify the kernel and image of
this map.

\begin{prop} Suppose $p>1$, $l\ge 0$ and\/ $\ga<0$ with\/
$[\ga,-\ga]\cap\D_{(\d+\d^*)_0}=[\ga,-\ga]\cap
\D_{(\d^*\d+\d\d^*)_0}=\{0\}$. Then the map
$\Ker\bigl((\d+\d^*)^p_{l+2,\ga}\bigr)\ra H^*(C,\R)$
given by $\chi\!\mapsto\![\chi]$ is injective, with image
that of the natural map~$H^*_{\rm cs}(C,\R)\!\ra\!H^*(C,\R)$.
\label{co3prop4}
\end{prop}

\begin{proof} Lockhart \cite[Ex.~0.14]{Lock} shows that
the vector space $\H^2(\La^kT^*C,g)$ of
closed, coclosed $k$--forms in $L^2(\La^kT^*C)$ on an
asymptotically cylindrical Riemannian manifold $(C,g)$
is isomorphic under $\chi\mapsto[\chi]$ with the image
of $H^k_{\rm cs}(C,\R)$ in $H^k(C,\R)$. Taking the
direct sum over $k=0,\ldots,n$, this implies that
for $l\ge 0$ the map
\begin{equation*}
\Ker\bigl((\d+\d^*)^2_{l+2,0}\bigr)\ra H^*(C,\R),
\qquad \chi\mapsto[\chi]
\end{equation*}
is injective, with image that of the natural
map~$H^*_{\rm cs}(C,\R)\ra H^*(C,\R)$.

Using Proposition \ref{co3prop1}, $\ga\notin\D_{(\d+\d^*)_0}$
and $\ga<0$ we have
\begin{equation*}
\Ker\bigl((\d+\d^*)^p_{l+2,\ga}\bigr)=
\Ker\bigl((\d+\d^*)^2_{l+2,\ga}\bigr)\subseteq
\Ker\bigl((\d+\d^*)^2_{l+2,0}\bigr).
\end{equation*}
Therefore $\chi\mapsto[\chi]$ is {\it injective} on
$\Ker\bigl((\d+\d^*)^p_{l+2,\ga}\bigr)$, with image
{\it contained in} that of $H^*_{\rm cs}(C,\R)\ra
H^*(C,\R)$. It remains to show $\chi\mapsto[\chi]$ is
{\it surjective} on this image.

Suppose $\eta\in H^j(C,\R)$ lies in the image of
$H^j_{\rm cs}(C,\R)$. Then we may write $\eta=[\phi]$
for $\phi$ a smooth, closed, compactly-supported $j$--form
on $C$. Hence $\d^*\phi\in\bigop_{k=0}^nL^p_{l+1,\ga}
(\La^kT^*C)$. We shall show that $\d^*\phi$ lies in the
image of \eq{co3eq11} with $l+1$ in place of $l$. Since
$\ga\notin\D_{(\d^*\d+\d\d^*)_0}$, as in Section~\ref{co31}
this holds if and only if $\langle\d^*\phi,\xi\rangle_{L^2}=0$
for all $\xi$ in~$\Ker\bigl((\d^*\d+\d\d^*)^q_{m+2,-\ga}\bigr)$.

But all such $\xi$ are {\it closed\/} by Proposition
\ref{co3prop3}, so $\langle\d^*\phi,\xi\rangle_{L^2}=
\langle\phi,\d\xi\rangle_{L^2}=0$. Therefore
$\d^*\phi=(\d^*\d+\d\d^*)\psi$ for some
$\psi\in L^p_{l+3,\ga}(\La^{j-1}T^*C)$. Hence
\begin{equation*}
(\d^*\d+\d\d^*)(\phi-\d\psi)=\d\bigl(\d^*\phi-
(\d^*\d+\d\d^*)\psi\bigr)=0,
\end{equation*}
and $\phi-\d\psi$ lies in $\Ker\bigl((\d^*\d+\d\d^*)^p_{
l+2,\ga}\bigr)$, which is $\Ker\bigl((\d+\d^*)^p_{l+2,\ga}
\bigr)$ by \eq{co3eq12}. As $[\phi-\d\psi]=[\phi]=\eta$
we have proved the surjectivity we need.
\end{proof}

\subsection{$\d_++\d^*$ on a 4--manifold}
\label{co33}

Now we restrict to $\dim C=4$, so that $(C,g)$ is an oriented
asymptotically cylindrical Riemannian 4--manifold asymptotic
to a Riemannian cylinder $(L\t\R,g_0)$. In Section~\ref{co4} we will
take $C$ to be an {\it asymptotically cylindrical coassociative
$4$--fold}. Consider the {\it asymptotically cylindrical linear
elliptic operator}
\begin{equation*}
\d_++\d^*\co C^\iy(\La^2_+T^*C)\op C^\iy(\La^4T^*C)\longra
C^\iy(\La^3T^*C).
\end{equation*}
Here $\d_+$ is the restriction of $\d$ to the self-dual
2--forms. We use this notation to distinguish $\d_++\d^*$
from $\d+\d^*$ in \eq{co3eq9}. Roughly speaking,
$\d_++\d^*$ is a {\it quarter} of $\d+\d^*$ in \eq{co3eq9},
as it acts on half of the even forms, rather than on all
forms. Its {\it formal adjoint\/} is
\begin{equation*}
\d^*_++\d\co C^\iy(\La^3T^*C)\longra C^\iy(\La^2_+T^*C)\op
C^\iy(\La^4T^*C),
\end{equation*}
where $\d^*_+$ is the projection of $\d^*$ to the self-dual
2--forms. We shall apply the results of Section~\ref{co32} to
study the extension
\begin{equation}
(\d_++\d^*)^p_{l+2,\ga}\co L^p_{l+2,\ga}(\La^2_+T^*C)\op
L^p_{l+2,\ga}(\La^4T^*C)\longra L^p_{l+1,\ga}(\La^3T^*C),
\label{co3eq17}
\end{equation}
for $p>1$, $l\ge 0$ and $\ga\in\R$. We begin with some
algebraic topology.

Suppose for simplicity that $C$ has no compact connected
components, so that $H^4(C,\R)=H^0_{\rm cs}(C,\R)=0$. Then
$L$ is a compact, oriented 3--manifold, and $C$ is the
interior $(\bar C)^\circ$ of a compact, oriented 4--manifold
$\bar C$ with boundary $\partial\bar C=L$. Thus we have a
long exact sequence in cohomology:
\begin{equation}
\begin{gathered}
\xymatrix@C14pt@R17pt{ 0 \ar[r] & H^0(C) \ar[r] & H^0(L) \ar[r]
&H^1_{\rm cs}(C) \ar[r]
& H^1(C) \ar[r]
& H^1(L) \ar[r]
& H^2_{\rm cs}(C) \ar[d] \\
0
& H^4_{\rm cs}(C) \ar[l]
& H^3(L) \ar[l]
& H^3(C) \ar[l]
&
H^3_{\rm cs}(C) \ar[l]
& H^2(L) \ar[l]
& H^2(C) \ar[l] }
\end{gathered}
\label{co3eq18}
\end{equation}
where $H^k(C)=H^k(C,\R)$ and $H^k(L)=H^k(L,\R)$ are
the de Rham cohomology groups, and $H^k_{\rm cs}(C,\R)$
is {\it compactly-supported de Rham cohomology}. Let
$b^k(C)$, $b^k(L)$ and $b^k_{\rm cs}(C)$ be the
corresponding Betti numbers.

By {\it Poincar\'e duality} we have
$H^k(C)\cong H^{4-k}_{\rm cs}(C)^*$ and $H^k(L)\cong
H^{3-k}(L)^*$, so that $b^k(C)=b^{4-k}_{\rm cs}(C)$ and
$b^k(L)=b^{3-k}(L)$. Note that \eq{co3eq18} is written
so that each vertically aligned pair of spaces are dual
vector spaces, and each vertically aligned pair of maps
are dual linear maps.

Let $V\subseteq H^2(C,\R)$ be the image of the
natural map $H^2_{\rm cs}(C,\R)\ra H^2(C,\R)$. Taking
alternating sums of dimensions in \eq{co3eq18} shows that
\begin{equation*}
\begin{split}
\dim V&=b^2_{\rm cs}(C)-b^1(L)+b^1(C)-b^1_{\rm cs}(C)-b^0(L)+b^0(C)\\
&=b^0(C)+b^1(C)+b^2(C)-b^3(C)-b^0(L)-b^1(L).
\end{split}
\end{equation*}
Now the cup product $\cup\co H^2_{\rm cs}(C,\R)\t H^2(C,\R)\ra\R$
restricted to $H^2_{\rm cs}(C,\R)\t V$ is zero on the product
of the kernel of $H^2_{\rm cs}(C,\R)\ra H^2(C,\R)$ with $V$.
Hence it pushes forward to a quadratic form $\cup\co V\t V\ra\R$,
which is {\it symmetric} and {\it nondegenerate}.

Suppose $V=V_+\op V_-$ is a decomposition of $V$ into subspaces
with $\cup$ positive definite on $V_+$ and negative definite
on $V_-$. Then $\dim V_+$ and $\dim V_-$ are {\it topological
invariants} of $C,L$. That is, they depend only on $C$ as an
oriented 4--manifold, and not on the choice of subspaces~$V_\pm$.

We now identify the kernel and cokernel of $(\d_++\d^*)^p_{l+2,\ga}$
in \eq{co3eq17} for small~$\ga<0$.

\begin{thm} Let\/ $(C,g)$ be an oriented, asymptotically
cylindrical Riemannian $4$--manifold asymptotic to $(L\t\R,g_0)$,
and use the notation of\/ Sections~\ref{co31}--\ref{co32}. Suppose
$\max\bigl(\D_{(\d_++\d^*)_0}\cap(-\iy,0)\bigr)<\ga<0$, and let\/
$p,q>1$ with\/ $\frac{1}{p}+\frac{1}{q}=1$ and\/ $l,m\ge 0$.
Then from Section~\ref{co31} the operator $(\d_++\d^*)^p_{l+2,\ga}$
of\/ \eq{co3eq17} is Fredholm with
\begin{equation}
\Coker\bigl((\d_++\d^*)^p_{l+2,\ga}\bigr)
\cong\Ker\bigl((\d^*_++\d)^q_{m+2,-\ga}\bigr)^*.
\label{co3eq19}
\end{equation}
The kernel $\Ker\bigl((\d_++\d^*)^p_{l+2,\ga}\bigr)$ is
a vector space of smooth, closed, self-dual\/ $2$--forms.
The map $\Ker\bigl((\d_++\d^*)^p_{l+2,\ga}\bigr)\ra H^2(C,\R)$,
$\chi\mapsto[\chi]$ induces an isomorphism of\/
$\Ker\bigl((\d_++\d^*)^p_{l+2,\ga}\bigr)$ with a maximal
subspace $V_+$ of the subspace $V\subseteq H^2(C,\R)$
defined above on which the cup product $\cup\co V\t V\ra\R$
is positive definite. Hence
\begin{equation}
\dim\Ker\bigl((\d_++\d^*)^p_{l+2,\ga}\bigr)=\dim V_+,
\label{co3eq20}
\end{equation}
which is a topological invariant of\/ $C,L$ from above.
Also, $\Ker\bigl((\d^*_++\d)^q_{m+2,-\ga}\bigr)$ is a
vector space of smooth, closed and coclosed\/ $3$--forms.
\label{co3thm2}
\end{thm}

\begin{proof} The first part follows immediately from
Section~\ref{co31}. As $\Ker\bigl((\d_++\d^*)^p_{l+2,\ga}\bigr)$
and $\Ker\bigl((\d^*_++\d)^q_{m+2,-\ga}\bigr)$ depend only on
the connected component of $\R\sm\D_{(\d_++\d^*)_0}$ containing
$\ga$, we can make $\md{\ga}$ smaller if necessary to ensure
that $[\ga,-\ga]\cap\D_{(\d+\d^*)_0}=[\ga,-\ga]\cap
\D_{(\d^*\d+\d\d^*)_0}=\{0\}$, using the notation of Section~\ref{co32}.

Suppose $(\ze^2_+,\ze^4)\in\Ker\bigl((\d_++\d^*)^p_{l+2,\ga}\bigr)$.
Then $\d\ze^2_++\d^*\ze^4=0$, so applying $*$ and noting that
$*\ze^2_+=\ze^2_+$ gives $\d^*\ze^2_+-\d(*\ze^4)=0$. Hence
$(\d+\d^*)(-*\ze^4+\ze^2_++\ze^4)=0$, that is, the mixed form
$-*\ze^4+\ze^2_++\ze^4$ lies in the kernel of \eq{co3eq10}.
Proposition \ref{co3prop3} now implies that $\ze^2_+$ and
$\ze^4$ are smooth, closed and coclosed, and therefore
$\ze^2_+$ and $\ze^4$ also lie in the kernel, and
$\Ker\bigl((\d_++\d^*)^p_{l+2,\ga}\bigr)\subseteq
\Ker\bigl((\d+\d^*)^p_{l+2,\ga}\bigr)$, where
$(\d+\d^*)^p_{l+2,\ga}$ is as in Section~\ref{co32}.

Since $H^4(C,\R)=0$, injectivity in Proposition \ref{co3prop4}
implies that $\ze^4=0$. Thus $\Ker\bigl((\d_++\d^*)^p_{l+2,\ga}
\bigr)$ is a vector space of smooth, closed, self-dual $2$--forms,
as we have to prove. Write $\H^2$ for the space of 2--forms
in $\Ker\bigl((\d+\d^*)^p_{l+2,\ga}\bigr)$. Then by Proposition
\ref{co3prop4} the map $\H^2\ra H^2(C,\R)$, $\chi\mapsto
[\chi]$ is injective with image that of $H^2_{\rm cs}(C,\R)$
in $H^2(C,\R)$, which is $V$ in the notation above.

Under this isomorphism, the cup product on $V$ is given by
\begin{equation}
\chi\cup\xi=\int_C\chi\w\xi=\int_C(\chi,*\xi)\d V_g=
\langle\chi,*\xi\rangle_{L^2},
\label{co3eq21}
\end{equation}
for $\chi,\xi\in\H^2$. The Hodge star $*$ maps
$\H^2\cong V$ to itself with $*^2=1$. Let $V_\pm$
be the $\pm 1$ eigenspaces of $*$ on $V$. Then $V=V_+\op V_-$,
and \eq{co3eq21} implies that $\cup$ is positive definite on
$V_+$ and negative definite on $V_-$. Hence $\dim V_+$ is a
{\it topological invariant}, from above. But $\Ker\bigl(
(\d_++\d^*)^p_{l+2,\ga}\bigr)$ is the self-dual 2--forms
in $\H^2$. Thus $\chi\mapsto[\chi]$ induces an
isomorphism of $\Ker\bigl((\d_++\d^*)^p_{l+2,\ga}\bigr)$
with $V_+$, as we have to prove.

Finally, suppose $\ze^3\in\Ker\bigl((\d^*_++\d)^q_{m+2,-\ga}
\bigr)$, so that $\d_+^*\ze^3=\d\ze^3=0$. Then
\begin{equation*}
0=\d^*(\d\ze^3)+2\d(\d_+^*\ze^3)=\d^*\d\ze^3
+\d(\d^*\ze^3-\d(*\ze^3))=(\d^*\d+\d\d^*)\ze_3,
\end{equation*}
so $\ze^3$ lies in $\Ker\bigl((\d^*\d+\d\d^*)^q_{m+2,-\ga}
\bigr)$, which consists of smooth closed and coclosed forms
by Proposition \ref{co3prop3}. Thus $\Ker\bigl((\d^*_++\d
)^q_{m+2,-\ga}\bigr)$ is a vector space of smooth, closed
and coclosed $3$--forms.
\end{proof}

We can say more about the cokernel \eq{co3eq19}. Its dimension is
\begin{equation*}
\dim\Ker\bigl((\d^*_++\d)^q_{m+2,-\ga}\bigr)=b^0(L)-b^0(C)+b^1(C).
\end{equation*}
The map $\Ker\bigl((\d^*_++\d)^q_{m+2,-\ga}\bigr)\ra H^3(C,\R)$,
$\chi\mapsto[\chi]$ is surjective, with kernel of dimension
$b^0(L)-b^0(C)+b^1(C)-b^3(C)\ge 0$, which is the dimension of
the kernel of $H^3_{\rm cs}(C,\R)\ra H^3(C,\R)$. So we can
think of $\Ker\bigl((\d^*_++\d)^q_{m+2,-\ga}\bigr)$ as a space
of closed and coclosed 3--forms filling out all of $H^3_{\rm cs}
(C,\R)$ and $H^3(C,\R)$. But we will not need these facts, so
we shall not prove them.

\subsection{Conditions on the rate $\ga$ for $\d_++\d^*$ to be Fredholm}
\label{co34}

Finally we determine the set $\D_{(\d_++\d^*)_0}$ for the
cylindrical operator $(\d_++\d^*)_0$ on $L\t\R$, which by
Theorem \ref{co3thm1} gives the set of $\ga$ for which
\eq{co3eq17} is not Fredholm.

\begin{prop} Let\/ $(C,g)$ be an oriented, asymptotically
cylindrical Riemannian $4$--manifold asymptotic to $(L\t\R,g_0)$,
let\/ $p>1$, $l\ge 0$ and\/ $\ga\in\R$, and define
$(\d_++\d^*)^p_{l+2,\ga}$ as in \eq{co3eq17}. Then
$(\d_++\d^*)^p_{l+2,\ga}$ is not Fredholm if and only if
either $\ga=0$, or $\ga^2$ is a positive eigenvalue of\/
$\De=\d^*\d$ on functions on $L$, or $\ga$ is an
eigenvalue of $-*\d$ on coexact\/ $1$--forms on~$L$.
\label{co3prop5}
\end{prop}

\begin{proof} Throughout the proof $*$, $\d^*$ mean the
Hodge star and $\d^*$ on $L$, not on $L\t\R$, and $\d^*_4$ is
$\d^*$ on $L\t\R$. An element of $C^\iy(\La^2_+T^*(L\t\R)\ot_\R\C)$
invariant under translations in $\R$ may be written uniquely
in the form $\chi\w\d t+*\chi$ for $\chi\in C^\iy(T^*L\ot_\R\C)$.
An element of $C^\iy(\La^4T^*(L\t\R)\ot_\R\C)$ invariant under
translations in $\R$ may be written uniquely as $f\,\d V_L\w\d t$
for $f\co L\ra\C$ smooth, where $\d V_L$ is the volume form on $L$.
By Definition \ref{co3def3} and Theorem \ref{co3thm1}, \eq{co3eq17}
is not Fredholm if and only if there exist $\de\in\R$ and $\chi,f$
as above and not both zero, satisfying
\begin{equation}
\d\bigl(e^{(\ga+i\de)t}(\chi\w\d t+*\chi)\bigr)+
\d^*_4\bigl(e^{(\ga+i\de)t}f\,\d V_L\w\d t\bigr)\equiv 0.
\label{co3eq22}
\end{equation}
Expanding \eq{co3eq22} yields
\begin{equation*}
e^{(\ga+i\de)t}\bigl[\d\chi\w\d t+(\ga+i\de)(*\chi)\w\d t+\d(*\chi)
-(*\d f)\w\d t+(\ga+i\de)f\,\d V_L\bigr]\equiv 0
\end{equation*}
on $L\t\R$, and separating components with and without $\d t$ gives
\begin{equation*}
\d\chi+(\ga+i\de)(*\chi)-(*\d f)\equiv 0
\quad\text{and}\quad
\d(*\chi)+(\ga+i\de)f\,\d V_L\equiv 0,
\end{equation*}
equations in 2-- and 3--forms on $L$ respectively. Applying the
Hodge star $*$ on $L$ shows that \eq{co3eq22} is equivalent
to the two equations
\begin{equation}
*\d\chi+(\ga+i\de)\chi-\d f\equiv 0
\quad\text{and}\quad
\d^*\chi-(\ga+i\de)f\equiv 0
\label{co3eq23}
\end{equation}
in 1--forms and functions on~$L$.

Since $\chi$ is a 1--form on $L$ which is a compact manifold,
one can use Hodge decomposition and write
$\chi=\chi_0\op\chi_1\op\chi_2$ where $\chi_0$ is
a harmonic 1--form, $\chi_1=\d f_1$ is an exact 1--form,
$f_1\in C^\iy(L)$ and $\chi_2=\d^*\eta$ is a co-exact
1--form, $\eta\in C^\iy(\La^2T^*L)$. Dividing the first equation
of \eq{co3eq23} into harmonic, exact and coexact components,
the system becomes
\begin{equation}
(\ga+i\de)\chi_0=0,\;\>
(\ga+i\de)\d f_1=\d f,\;\>
-*\d\chi_2=(\ga+i\de)\chi_2,\;\>
\d^*\d f_1=(\ga+i\de)f.
\label{co3eq24}
\end{equation}
The second and fourth equations of \eq{co3eq24} give
$\d^*\d f=(\ga+i\de)^2f$, and substituting the third
equation into itself gives $\d^*\d\chi_2=*\d*\d\chi_2=
(\ga+i\de)^2\chi_2$. But $\d^*\chi_2=0$ by definition. Thus
$\chi_0,\chi_1,\chi_2,f$ satisfy the equations
\begin{equation}
\begin{aligned}
(\ga+i\de)\chi_0&=0,&\quad
(\ga+i\de)\chi_1&=\d f,\\
(\d^*\d+\d\d^*)\chi_2&=(\ga+i\de)^2\chi_2,&\quad
\d^*\d f&=(\ga+i\de)^2f.
\end{aligned}
\label{co3eq25}
\end{equation}
When $\ga=\de=0$, $\chi=0$ and $f\equiv 1$ are a solution,
so \eq{co3eq17} is not Fredholm. So suppose $\ga+i\de\ne 0$.
Then $\chi_0=0$, and $\chi_1=0$ if $f=0$ by \eq{co3eq25},
so as $\chi,f$ are not both zero either $f\ne 0$ or
$\chi_2\ne 0$. If $f\ne 0$ then $(\ga+i\de)^2$ is a nonzero
eigenvalue of $\De=\d^*\d$ on functions by the last equation
of \eq{co3eq25}, so $\de=0$ as such eigenvalues are positive,
and $\ga^2$ is a positive eigenvalue of $\De$ as we want.
Conversely, if $\De f=\ga^2f$ for nonzero $\ga,f$ then
$\chi=\ga^{-1}\d f$ satisfies the equations, so \eq{co3eq17}
is not Fredholm.

If $\chi_2\ne 0$ then $(\ga+i\de)^2$ is a nonzero eigenvalue
of $\De=\d^*\d+\d\d^*$ on coexact 1--forms by the third equation
of \eq{co3eq25}, so $\de=0$ as above. The third equation of
\eq{co3eq24} then shows that $\ga$ is an eigenvalue of
$-*d$ on coexact 1--forms on $L$. Conversely, taking $\chi_2$
to be an eigenvector of $-*d$ on coexact 1--forms with eigenvalue
$\ga$ and $\chi_0=\chi_1=f=0$ solves the equations, so \eq{co3eq17}
is not Fredholm.
\end{proof}

\section{Proof of Theorem \ref{co1thm}}
\label{co4}

We now prove Theorem \ref{co1thm}. Let $(M,\vp,g)$ be an
asymptotically cylindrical $G_2$--manifold asymptotic to
$X\t(R,\iy )$, $R>0$, with decay rate $\al<0$. Let $C$ be
an asymptotically cylindrical coassociative 4--fold in $X$
asymptotic to $L\t(R',\iy )$ for $R'>R$ with decay rate
$\be$ for $\al\le\be<0$. Write $g_C=g\vert_C$ for the
metric on $C$, $g_X$ for the Calabi--Yau metric on $X$,
and $g_L=g_X\vert_L$ for the metric on $L$. Then
$(C,g_C)$ is an {\it asymptotically cylindrical
Riemannian $4$--manifold}, with {\it rate}~$\be$.

Suppose $\ga$ satisfies $\be<\ga<0$, and $(0,\ga^2]$
contains no eigenvalues of the Laplacian $\De_L$ on
functions on $L$, and $[\ga,0)$ contains no eigenvalues
of the operator $-*\d$ on coexact 1--forms on $L$. Let
$p>4$ and $l\ge 1$, and define
\begin{equation}
(\d_++\d^*)^p_{l+2,\ga}\co
L^p_{l+2,\ga}(\La^2_+T^*C)\op L^p_{l+2,\ga}(\La^4T^*C)
\longra L^p_{l+1,\ga}(\La^3T^*C)
\label{co4eq1}
\end{equation}
as in Section~\ref{co33}. Then Proposition \ref{co3prop5} and
the conditions on $\ga$ imply that $[\ga,0)\cap\D_{(\d_+
+\d^*)_0}=\emptyset$. Hence $\ga\notin\D_{(\d_++\d^*)_0}$,
so that $(\d_++\d^*)^p_{l+2,\ga}$ is {\it Fredholm} by
Theorem \ref{co3thm1}. Also, Theorem \ref{co3thm2} applies
to~$(\d_++\d^*)^p_{l+2,\ga}$.

Let $\nu_L$ be the normal bundle of $L$ in $X$, regarded as the
{\it orthogonal subbundle} to $TL$ in $TX\vert_L$, and
$\exp_L\co \nu_L\ra X$ the exponential map. For $r>0$, write
$B_r(\nu_L)$ for the subbundle of $\nu_L$ with fibre at $x$ the
open ball about 0 in $\nu_L\vert_x$ with radius $r$. Then for
small $\ep>0$, there is a {\it tubular neighbourhood\/} $T_L$ of
$L$ in $X$ such that $\exp_L\co B_{2\ep}(\nu_L)\ra T_L$ is a
diffeomorphism. Also, $\nu_L\t\R\ra L\t\R$ is the normal bundle to
$L\t\R$ in $X\t\R$ with exponential map $\exp_L\t{\mathop{\rm
id}}\co \nu_L\t\R\ra X\t\R$. Then $T_L\t\R$ is a tubular neighborhood
of $L\t\R$ in $X\t\R$, and $\exp_L\t{\mathop{\rm id}}\co B_{2\ep}
(\nu_L)\t\R\ra T_L\t\R$ is a diffeomorphism.

Let $K,R$, $\Psi\co X\t(R,\iy)\ra M\sm K$, and $K',R'>R$,
$\Phi\co L\t(R',\iy)\ra C\sm K'$, and the normal vector
field $v$ on $L\t(R',\iy)$ be as in Section~\ref{co23}, so
that \eq{co2eq3} commutes. Then $v$ is a section of
$\nu_L\t(R',\iy)\ra L\t(R',\iy)$, decaying at rate
$O(e^{\be t})$. Therefore making $K'$ and $R'$ larger
if necessary, we can suppose the graph of $v$ lies
in~$B_\ep(\nu_L)\t(R',\iy)$.

Write $\pi\co B_\ep(\nu_L)\t(R',\iy)\ra L\t(R',\iy)$ for
the natural projection. Define
\begin{equation}
\Xi\co B_\ep(\nu_L)\t(R',\iy)\ra M \quad\text{by}\quad
\Xi\co w\mapsto \Psi\bigl[(\exp_L\t{\mathop{\rm id}})
(v\vert_{\pi(w)}+w)\bigr].
\label{co4eq2}
\end{equation}
Here $w$ is a point in $B_\ep(\nu_L)\t(R',\iy)$,
in the fibre over $\pi(w)\in L\t(R',\iy)$. Thus
$v\vert_{\pi(w)}$ is a point in the same fibre,
which is a ball of radius $\ep$ in a vector space.

Hence $v\vert_{\pi(w)}+w$ lies in the open ball of
radius $2\ep$ in the same vector space, that is, in
the fibre of $B_{2\ep}(\nu_L)\t(R',\iy)$ over $\pi(w)$.
Therefore $(\exp_L\t{\mathop{\rm id}})(v\vert_{\pi(w)}+w)$
is well-defined and lies in $T_L\t(R',\iy)\subset
X\t(R',\iy)$, and $\Xi(w)$ is well-defined. Since
$\exp_L\t{\mathop{\rm id}}\co B_{2\ep}(\nu_L)\t\R\ra
T_L\t\R$ and $\Psi\co X\t(R,\iy)\ra M\sm K$ are
diffeomorphisms, we see that $\Xi$ is a diffeomorphism
with its image.

Identify $L\t(R',\iy)$ with the zero section in
$B_\ep(\nu_L)\t(R',\iy)$. Then as \eq{co2eq3} commutes
we see that $\Xi\vert_{L\t(R',\iy)}\equiv\Phi$.
Let $\nu_C$ be the normal bundle of $C$ in $M$, which
we regard not as the orthogonal subbundle to $TC$ in
$TM\vert_C$, but rather as the {\it quotient bundle}
$TM\vert_C/TC$. Define an isomorphism $\xi$ between the
vector bundles $\nu_L\t(R',\iy)$ and $\Phi^*(\nu_C)$
over $L\t(R',\iy)$ as follows.

As $\Xi$ is a diffeomorphism with its image,
$\d\Xi\co T\bigl(B_\ep(\nu_L)\t(R',\iy)\bigr)\ra\Xi^*(TM)$
is an isomorphism. Restricting this to the zero section
$L\t(R',\iy)$ and noting that $\Xi\vert_{L\t(R',\iy)}
\equiv\Phi$, we see that
\begin{equation}
\d\Xi\vert_{L\t(R',\iy)}\co T\bigl(B_\ep(\nu_L)\t(R',\iy)\bigr)
\vert_{L\t(R',\iy)}\longra\Phi^*(TM)
\label{co4eq3}
\end{equation}
is an isomorphism. As $\Phi\co L\t(R',\iy)\ra C\sm K'$ is a
diffeomorphism,
\begin{equation}
\d\Phi\co T\bigl(L\t(R',\iy)\bigr)\longra\Phi^*(TN)
\label{co4eq4}
\end{equation}
is an isomorphism. But \eq{co4eq4} is the restriction
of \eq{co4eq3} to a vector subbundle. Quotienting
\eq{co4eq3} by \eq{co4eq4} gives an isomorphism
\begin{multline}
\xi=\d\Xi\vert_{L\t(R',\iy)}\co \nu_L\t(R',\iy)\cong
\frac{T\bigl(B_\ep(\nu_L)\t(R',\iy)\bigr)\vert_{L\t(R',\iy)}
}{T\bigl(L\t(R',\iy)\bigr)}\\
\longra\Phi^*(TM)/\Phi^*(TN)\cong\Phi^*(TM/TN)=\Phi^*(\nu_C).
\label{co4eq5}
\end{multline}
Now choose a small $\ep'>0$, a tubular neighborhood $T_C$ of $C$
in $M$, and a diffeomorphism $\Th\co B_{\ep'}(\nu_C)\ra T_C$
satisfying the conditions:
\begin{itemize}
\item[(i)] $\Th\vert_C\equiv{\mathop{\rm id}}_C\co C\ra C$, where
$C\subset B_{\ep'}(\nu_C)$ is the zero section.
\item[(ii)] By (i), $\d\Th\vert_C\co T(B_{\ep'}(\nu_C))\vert_C\ra
TM\vert_C$ is an isomorphism, which restricts to the identity
on the subbundles $TC$ of each side. Hence it induces an
isomorphism $\d\Th\vert_C\co \nu_C\cong T\bigl(B_\ep(\nu_C)
\bigr)\vert_C/TC\ra TM\vert_C/TC=\nu_C$.
This isomorphism is the identity map.
\item[(iii)] $\ep'$ is small enough that $\xi^*\bigl(B_{\ep'}
(\nu_C)\bigr)\subset B_\ep(\nu_L)\t(R',\iy)\subset \nu_L\t
(R',\iy)$, and $\Th\ci\xi\equiv\Xi$ on~$\xi^*\bigl(B_{\ep'}
(\nu_C)\bigr)$.
\end{itemize}
Notice that (iii) determines $\Th$ and $T_C$ uniquely on
$B_{\ep'}(\nu_C)\vert_{C\sm K'}$, and by construction
here it satisfies (i) and (ii). Thus, it remains only to
choose $T_C$ and $\Th$ satisfying (i), (ii) over the
compact set $K'\subset C$, which is possible by standard
differential topology.

The point of all this is that we have chosen a local
identification $\Th$ between $\nu_C$ and $M$ near $C$
that is compatible in a nice way with the asymptotic
identifications $\Phi,\Psi$ of $C$, $M$ with $L\t\R$
and $X\t\R$. Using $\Th$, submanifolds $\ti C$ of
$M$ close to $C$ are identified with small sections
$s$ of $\nu_C$, and importantly, the asymptotic
convergence of $\ti C$ to $C$, and so to $L\t\R$, is
reflected in the asymptotic convergence of $s$ to~0.

As in the proof of Theorem \ref{co2thm}, the map
$V\mapsto(V\cdot\vp\vert_x)\vert_{T_xC}$ defines an
{\it isomorphism} $\nu_C\ra\La^2_+T^*C$. (Note that
since $\vp\vert_C\equiv 0$, this map is well-defined
for $V\in T_xM/T_xC$, rather than just for $V\in T_xM$
orthogonal to $T_xC$.) We now {\it identify} $\nu_C$
with $\La^2_+T^*C$, and regard $\Th$ as a map
$\Th\co B_{\ep'}(\La^2_+T^*C)\ra T_C\subset M$.

Write $L^p_{l+2,\ga}\bigl(B_{\ep'}(\La^2_+T^*C)\bigr)$
for the subset of $\ze^2_+\in L^p_{l+2,\ga}(\La^2_+T^*C)$
which are sections of $B_{\ep'}(\La^2_+T^*C)$, that is,
$\md{\ze^2_+}<\ep'$ on $C$. Since $L^p_{l+1,\ga}\hookrightarrow
C^0$ by Sobolev embedding this is an open condition on $\ze^2_+$,
so $L^p_{l+2,\ga}\bigl(B_{\ep'}(\La^2_+T^*C)\bigr)$ is an
open subset of the Banach space~$L^p_{l+2,\ga}(\La^2_+T^*C)$.

Define $Q\co L^p_{l+2,\ga}\bigl(B_{\ep'}(\La^2_+T^*C)\bigr)\ra
\{$3--forms on $C\}$ by $Q(\ze^2_+)=(\Th\ci\ze^2_+)^*(\vp)$.
That is, we regard the section $\ze^2_+$ as a map $C\ra B_{\ep'}
(\La^2_+T^*C)$, so $\Th\ci\ze^2_+$ is a map $C\ra T_C\subset M$,
and thus $(\Th\ci\ze^2_+)^*(\vp)$ is a 3--form on $C$.
The point of this definition is that if
$\Ga_{\ze^2_+}$ is the {\it graph} of $\ze^2_+$ in
$B_{\ep'}(\La^2_+T^*C)$ and $\ti C=\Th(\Ga_{\ze^2_+})$
its image in $M$, then $\ti C$ is {\it coassociative}
if and only if $\vp\vert_{\ti C}\equiv 0$, which holds
if and only if $Q(\ze^2_+)=0$. So $Q^{-1}(0)$ parametrizes
coassociative 4--folds $\ti C$ close to~$C$.

We now consider which class of 3--forms $Q$ maps to.

\begin{prop} $Q\co L^p_{l+2,\ga}\bigl(B_{\ep'}(\La^2_+T^*C)
\bigr)\longra L^p_{l+1,\ga}(\La^3T^*C)$ is a smooth map
of Banach manifolds. The linearization of\/ $Q$ at\/ $0$
is\/~$\d Q(0)\co \ze^2_+\mapsto\d\ze^2_+$.
\label{co4prop1}
\end{prop}

\begin{proof} As in the proof of Theorem \ref{co2thm}, the
functional form of $Q$ is
\begin{equation*}
Q(\ze^2_+)\vert_x=F\bigl(x,\ze^2_+\vert_x,\na\ze^2_+\vert_x\bigr)
\quad\text{for $x\in C$,}
\end{equation*}
where $F$ is a smooth function of its arguments. Since $p>4$
and $l\ge 1$ we have $L^p_{l+2,\ga}(\La^2_+T^*C)\hookrightarrow
C^1_\ga(\La^2_+T^*C)$ by Sobolev embedding. General arguments
then show that {\it locally} $Q(\ze^2_+)$ is $L^p_{l+1}$.

To show $Q(\ze^2_+)$ lies in $L^p_{l+1,\ga}(\La^3T^*C)$, we must
know something of the asymptotic behavior of $F$ at infinity.
Essentially $Q(\ze^2_+)$ is the restriction to $\Ga(\ze^2_+)$ of
the 3--form $\Th^*(\vp)$ on $B_{\ep'} (\La^2_+T^*C)$. Using the
identifications $\Phi\co L'\t(R'\iy)\ra C\sm K'$ and
$\xi\co \nu_L\t(R',\iy) \ra\Phi^*(\La^2_+T^*C)$ over $C\sm K'$, by
(iii) above this 3--form becomes $\Xi^*(\vp)$ on~$B_\ep(\nu_L)\t
(R',\iy)$.

But the asymptotic conditions on $\Phi,\Psi$ and $v$ imply that
$\Xi^*(\vp)$ is the sum of a translation-invariant 3--form on
$B_\ep(\nu_L)\t\R$, the pullback of the cylindrical $G_2$ form
$\vp_0$ on $X\t\R$, and an error term which decays at rate
$O(e^{\be t})$, with all its derivatives. As $\be<\ga<0$, it is
not difficult to see from this that $Q$ maps to
$L^p_{l+1,\ga}(\La^3T^*C)$. Smoothness of $Q$ holds by general
principles. Finally, the linearization of $Q$ is $\d$, by the
calculation of McLean alluded to in Theorem \ref{co2thm}. As the
calculation is local, it does not matter that we are on a
noncompact manifold~$C$.
\end{proof}

Next we show that the image of $Q$ consists of {\it exact\/}
3--forms. Notice that we lose one degree of differentiability:
although the image of $Q$ lies in $L^p_{l+1,\ga}(\La^3T^*C)$,
we claim only that it lies in the exact 3--forms
in~$L^p_{l,\ga}(\La^3T^*C)$.

\begin{prop}
$$Q\bigl(L^p_{l+2,\ga}(B_{\ep'}(\La^2_+T^*C))
\bigr)\!\subseteq\!\d\bigl(L^p_{l+1,\ga}(\La^2T^*C)\bigr)
\!\subset\!L^p_{l,\ga}(\La^3T^*C).\vadjust{\vskip-20pt}$$
\label{co4prop2}
\end{prop}
\begin{proof} Consider the restriction of the 3--form
$\vp$ to the tubular neighborhood $T_C$ of $C$. As $\vp$ is
closed, and $T_C$ retracts onto $C$, and $\vp\vert_C\equiv 0$, we
see that $\vp\vert_{T_C}$ is {\it exact}. Thus we may write
$\vp\vert_{T_C}= \d\th$ for $\th\in C^\iy(\La^2T^*T_C)$. Since
$\vp\vert_C\equiv 0$ we may choose $\th\vert_C\equiv 0$. Also, as
$\vp$ is asymptotic to $O(e^{\be t})$ with all its derivatives to
a translation-invariant 3--form $\vp_0$ on $X\t\R$, we may take
$\th$ to be asymptotic to $O(e^{\be t})$ with all its derivatives
to a translation-invariant 2--form on~$T_L\t\R$.

The proof of Proposition \ref{co4prop1} now shows that
the map $\ze^2_+\mapsto(\Th\ci\ze^2_+)^*(\th)$ maps
$L^p_{l+2,\ga}\bigl(B_{\ep'}(\La^2_+T^*C)\bigr)\ra
L^p_{l+1,\ga}(\La^2T^*C)$. But
\begin{equation*}
Q(\ze^2_+)=(\Th\ci\ze^2_+)^*(\vp)=(\Th\ci\ze^2_+)^*(\d\th)=
\d\bigl[(\Th\ci\ze^2_+)^*(\th)\bigr],
\end{equation*}
so $Q(\ze^2_+)\in\d\bigl(L^p_{l+1,\ga}(\La^2T^*C)\bigr)$
for~$\ze^2_+\in L^p_{l+2,\ga}(B_{\ep'}(\La^2_+T^*C))$.
\end{proof}

As in the proof of Theorem \ref{co2thm}, we augment $Q$ by a
space of 4--forms on $C$ to make it {\it elliptic}. Define
\begin{gather*}
P\co L^p_{l+2,\ga}\bigl(B_{\ep'}(\La^2_+T^*C)
\bigr)\t L^p_{l+2,\ga}(\La^4T^*C)\longra
L^p_{l+1,\ga}(\La^3T^*C)\\
\text{by}\qquad
P(\ze^2_+,\ze^4)=Q(\ze^2_+)+\d^*\ze^4.
\end{gather*}
Proposition \ref{co4prop1} implies that the linearization
$\d P(0,0)$ of $P$ at 0 is the Fredholm operator
$(\d_++\d^*)^p_{l+2,\ga}$ of \eq{co4eq1}. Define
$\mathcal C$ to be the image of $(\d_++\d^*)^p_{l+2,\ga}$.
Then $\mathcal C$ is a {\it Banach subspace} of
$L^p_{l+1,\ga}(\La^3T^*C)$, since $(\d_++\d^*)^p_{l+2,\ga}$
is Fredholm. We show $P$ maps into~$\mathcal C$.

\begin{prop} $P$ maps $L^p_{l+2,\ga}\bigl(B_{\ep'}(\La^2_+T^*C)
\bigr)\t L^p_{l+2,\ga}(\La^4T^*C)\longra\mathcal C$.
\label{co4prop3}
\end{prop}

\begin{proof} Let $(\ze^2_+,\ze^4)\in L^p_{l+2,\ga}
\bigl(B_{\ep'}(\La^2_+T^*C)\bigr)\t L^p_{l+2,\ga}
(\La^4T^*C)$ so that $P(\ze^2_+,\ze^4)\allowbreak=Q(\ze^2_+)+\d^*\ze^4$
lies in $L^p_{l+1,\ga}(\La^3T^*C)$. We must show it lies
in $\mathcal C$. Since $\ga\notin\D_{(\d_++\d^*)_0}$, from
Section~\ref{co3} this holds if and only if
\begin{equation}
\big\langle Q(\ze^2_+)+\d^*\ze^4,\chi\big\rangle_{L^2}=0
\qquad\text{for all $\chi\in\Ker\bigl((\d^*_++\d)^q_{
m+2,-\ga}\bigr)$,}
\label{co4eq6}
\end{equation}
where $\frac{1}{p}+\frac{1}{q}=1$ and~$m\ge 0$.

By Theorem \ref{co3thm2}, $\Ker\bigl((\d^*_++\d)^q_{
m+2,-\ga}\bigr)$ consists of {\it closed and coclosed\/}
3--forms $\chi$. By Proposition \ref{co4prop2} we have
$Q(\ze^2_+)=\d\la$ for $\la\in L^p_{l+1,\ga}(\La^2T^*C)$. So
\begin{equation*}
\big\langle Q(\ze^2_+)+\d^*\ze^4,\chi\big\rangle_{L^2}=
\big\langle \d\la,\chi\big\rangle_{L^2}+
\big\langle \d^*\ze^4,\chi\big\rangle_{L^2}=
\big\langle \la,\d^*\chi\big\rangle_{L^2}+
\big\langle \ze^4,\d\chi\big\rangle_{L^2}=0
\end{equation*}
for $\chi\in\Ker\bigl((\d^*_++\d)^q_{m+2,-\ga}\bigr)$,
as $\chi$ is closed and coclosed, and the inner products
and integration by parts are valid because of the matching
of rates $\ga,-\ga$ and $L^p,L^q$ with $\frac{1}{p}+
\frac{1}{q}=1$. So \eq{co4eq6} holds, and $P$ maps
into~$\mathcal C$.
\end{proof}

We now apply the {\it Implicit Mapping Theorem for Banach
spaces}, \cite[Theorem~1.2.5]{Joyc1}.

\begin{thm} Let\/ $\mathcal{A},\mathcal{B}$ and\/ $\mathcal{C}$ be
Banach spaces, and\/ $\mathcal{U},\mathcal{V}$ open neighborhoods
of\/ $0$ in $\mathcal{A}$ and\/ $\mathcal{B}$. Suppose that the
function $P\co \mathcal{U}\t\mathcal{V}\ra\mathcal{C}$ is smooth
with\/ $P(0,0)=0$, and that\/ $\d P_{(0,0)}|_\mathcal{B}\co
\mathcal{B}\ra\mathcal{C}$ is an isomorphism of\/
$\mathcal{B},\mathcal{C}$ as vector and topological spaces. Then
there exists a connected open neighbourhood\/ $\mathcal{U}'
\subset\mathcal{U}$ of\/ $0$ in $\mathcal{A}$ and a unique smooth
map $G\co \mathcal{U}'\ra\mathcal{V}$ such that\/ $G(0)=0$ and\/
$P(x,G(x))=0$ for all\/~$x\in\mathcal{U}'$. \label{co4thm}
\end{thm}

Define $\mathcal{A}$ to be
$\Ker\bigl((\d_++\d^*)^p_{l+2,\ga}\bigr)$, and $\mathcal{B}$ to be
the subspace of $L^p_{l+2,\ga}(\La^2_+T^*C)\allowbreak \op
L^p_{l+2,\ga}(\La^4T^*C)$ which is $L^2$--orthogonal to $\mathcal{A}$. As
$\mathcal A$ is finite-dimensional and the $L^2$ inner product is
continuous on $L^p_{l+2,\ga}(\La^2_+T^*C)\op L^p_{l+2,\ga}
(\La^4T^*C)$, both $\mathcal{A},\mathcal{B}$ are {\it Banach
spaces}, and $\mathcal{A}\op\mathcal{B}=L^p_{l+2,\ga}
(\La^2_+T^*C)\op L^p_{l+2,\ga}(\La^4T^*C)$. Choose open
neighborhoods $\mathcal{U},\mathcal{V}$ of 0 in $\mathcal{A},
\mathcal{B}$ such that $\mathcal{U}\t\mathcal{V}\subseteq
L^p_{l+2,\ga}\bigl(B_{\ep'}(\La^2_+T^*C)\bigr)\t L^p_{l+2,
\ga}(\La^4T^*C)$.

Let $P,\mathcal{C}$ be as above. Then by Propositions
\ref{co4prop1} and \ref{co4prop3}, $P\co \mathcal{U}\t
\mathcal{V}\ra\mathcal{C}$ is a smooth map of Banach
manifolds with $P(0,0)=0$, and linearization $\d P(0,0)=
(\d_++\d^*)^p_{l+2,\ga}\co \mathcal{A}\op\mathcal{B}\ra
\mathcal{C}$. By definition $\mathcal{A}$ is the kernel
and $\mathcal{C}$ is the image of $(\d_++\d^*)^p_{l+2,\ga}$.
Hence $\d P_{(0,0)}|_\mathcal{B}\co \mathcal{B}\ra\mathcal{C}$
is an isomorphism of vector spaces, and as it is a
continuous linear map of Banach spaces, it is an isomorphism
of topological spaces by the Open Mapping Theorem.

Thus Theorem \ref{co4thm} applies, and gives a connected open
neighborhood $\mathcal{U}'$ of 0 in $\mathcal{U}$, and a smooth
map $G\co \mathcal{U}'\ra\mathcal{V}$ such that $G(0)=0$ and
$P(x,G(x))\equiv 0$. Moreover $P^{-1}(0)$ coincides with
$\bigl\{(x,G(x)):x\in\mathcal{U}'\bigr\}$ near $(0,0)$, and so is
smooth, finite-dimensional and locally isomorphic
to~$\mathcal{A}=\Ker\bigl((\d_++\d^* )^p_{l+2,\ga}\bigr)$.

\begin{lem} $P^{-1}(0)=Q^{-1}(0)\t\{0\}$.
\label{co4lem1}
\end{lem}

\begin{proof} Clearly $Q^{-1}(0)\t\{0\}\subseteq P^{-1}(0)$.
Suppose $(\ze^2_+,\ze^4)\in P^{-1}(0)$, so that $Q(\ze^2_+)
+\d^*\ze^4=0$. We shall show that $\ze^4=0$, so that
$Q(\ze^2_+)=0$, and thus $P^{-1}(0)\subseteq Q^{-1}(0)\t\{0\}$.
By Proposition \ref{co4prop2} we have $Q(\ze^2_+)=\d\la$ for
$\la\in L^p_{l+1,\ga}(\La^2T^*C)$, so $\d\la=-\d^*\ze^4$. Hence
\begin{equation*}
\lnm{\d^*\ze^4}{2}^2=
\big\langle \d^*\ze^4,\d^*\ze^4\big\rangle_{L^2}=
-\big\langle \d^*\ze^4,\d\la\big\rangle_{L^2}=
-\big\langle\ze^4,\d^2\la\big\rangle_{L^2}=0,
\end{equation*}
where the inner products and integration by parts are
valid as $L^p_{l+2,\ga}\hookrightarrow L^2_2$. Thus
$Q(\ze^2_+)=\d^*\ze^4=0$. But $\d^*\ze^4\cong\na\ze^4$ as
$\ze^4$ is a 4--form, so $\ze^4$ is {\it constant}. Since also
$\ze^4\ra 0$ near infinity in $C$, we have~$\ze^4\equiv 0$.
\end{proof}

Combining this with the previous description of $P^{-1}(0)$
shows that $Q^{-1}(0)$ is smooth, finite-dimensional and
locally isomorphic to $\Ker\bigl((\d_++\d^*)^p_{l+2,\ga}\bigr)$.
Next we show that $Q^{-1}(0)$ is independent of $l$, and so
consists of smooth solutions.

\begin{prop} If\/ $\ze^2_+\in Q^{-1}(0)$ then $\ze^2_+\in
L^p_{m+2,\ga}(\La^2_+T^*C)$ for all\/~$m\ge 1$.
\label{co4prop4}
\end{prop}

\begin{proof} Let $m\ge 1$ and $(\ze^2_+,\ze^4)\in L^p_{m+2,\ga}
\bigl(B_{\ep'}(\La^2_+T^*C)\bigr)\t L^p_{m+2,\ga}(\La^4T^*C)$,
so that $P(\ze^2_+,\ze^4)=Q(\ze^2_+)+\d^*\ze^4$ lies in
$L^p_{m+1,\ga}(\La^3T^*C)$. Write $\na$ for the Levi--Civita
connection of $g_C$ on $C$. We shall apply the $\na$--Laplacian
$\De_C=g_C^{ab}\na_a\na_b$ to $P(\ze^2_+,\ze^4)$, using the
index notation. Writing $Q$ in terms of $F$ as in Proposition
\ref{co4prop1}, we have
\begin{equation}
\De_C\bigl[P(\ze^2_+,\ze^4)\bigr]=g_C^{ab}\na_a\na_b
\bigl[F\bigl(x,\ze^2_+\vert_x,\na\ze^2_+\vert_x\bigr)\bigr]
+\De_C\d^*\ze^4.
\label{co4eq7}
\end{equation}
We shall expand \eq{co4eq7} in terms of the derivatives of $F$ by
the chain rule, using the following notation. Write $F=F(x,y,z)$.
Write $\na^xF$ for the derivative of $F$ `in the $x$ direction'
using $\na$. That is, $(\na^xF)(x',y',z')=\na
\bigl(F(x,y(x),z(x)\bigr)\vert_{x'}$, where $y(x),z(x)$ satisfy
$y(x')=y'$, $z(x')=z'$ and $(\na y)(x')=(\na z)(x')=0$. Write
$\pd^yF$ and $\pd^zF$ for the partial derivatives of $F$ in the
$y,z$ directions. That is, with $x$ {\it held constant\/} the
domains $\La^2_+T^*_xC$ and $T^*_xC\ot\La^2_+T^*_xC$ of $y,z$ and
the range $\La^3T^*_xC$ of $F$ are {\it vector spaces}, so the map
$(y,z)\mapsto F(x,y,z)$ is a smooth map between vector spaces, and
has well-defined partial derivatives.

Then expanding \eq{co4eq7} gives
\begin{equation}
\De_C\bigl[P(\ze^2_+,\ze^4)\bigr]=L_{\ze^2_+}(\ze^2_+,\ze^4)+
E(x,\ze^2_+,\na\ze^2_+,\na^2\ze^2_+),
\label{co4eq8}
\end{equation}
where $L_{\ze^2_+}$ is the linear third-order operator
for $0\le n\le m+1$ given by
\begin{equation}
\begin{gathered}
L_{\ze^2_+}\co L^p_{n+3,\ga}(\La^2_+T^*C\op\La^4T^*C)
\longra L^p_{n,\ga}(\La^3T^*C),\\
L_{\ze^2_+}(\xi^2_+,\xi^4)\vert_x=(\pd^zF)(x,\ze^2_+\vert_x,
\na\ze^2_+\vert_x)\cdot\De_C\na\xi^2_+\vert_x+\De_C\d^*\xi^4\vert_x,
\end{gathered}
\label{co4eq9}
\end{equation}
and $E(x,\ze^2_+,\na\ze^2_+,\na^2\ze^2_+)$ is the `error
term' given by
\begin{equation}
\begin{gathered}
E(x,\ze^2_+,\na\ze^2_+,\na^2\ze^2_+)=(\pd^yF)(x,y,z)\cdot\De_C\ze^2_+
+g_C^{ab}(\na^x_a\na^x_bF)(x,y,z)
\\
+g_C^{ab}(\na^x_a\pd^yF)(x,y,z)\cdot\na_b\ze^2_+
+g_C^{ab}(\na^x_a\pd^zF)(x,y,z)\cdot\na_b\na\ze^2_+
\\
+g_C^{ab}(\pd^y\na^x_bF)(x,y,z)\cdot\na_a\ze^2_+
+g_C^{ab}(\pd^y\pd^yF)(x,y,z)\cdot(\na_a\ze^2_+\ot\na_b\ze^2_+)
\\
+g_C^{ab}(\pd^y\pd^zF)(x,y,z)\cdot(\na_a\ze^2_+\ot\na_b\na\ze^2_+)
+g_C^{ab}(\pd^z\na^x_bF)(x,y,z)\cdot\na_a\na\ze^2_+
\\
+\!g_C^{ab}(\pd^z\pd^yF)(x,y,z){\cdot}(\na_a\na\ze^2_+{\ot}\na_b
\ze^2_+) \\
{+}g_C^{ab}(\pd^z\pd^zF)(x,y,z){\cdot}
(\na_a\na\ze^2_+{\ot}\na_b\!\na\ze^2_+).
\end{gathered}
\label{co4eq10}
\end{equation}
Here in \eq{co4eq9} and \eq{co4eq10} we have used `$\,\cdot\,$'
to denote various natural bilinear products, and in \eq{co4eq10}
we write $y=\ze^2_+\vert_x$ and $z=\na\ze^2_+\vert_x$ in the
arguments of~$F$.

The point is that \eq{co4eq8} splits $\De_C[P(\ze^2_+,
\ze^4)]$ up into the piece $L_{\ze^2_+}(\ze^2_+,\ze^4)$ containing
all the {\it third derivative} terms of $\ze^2_+,\ze^4$, plus
a term $E(x,\ze^2_+,\na\ze^2_+,\na^2\ze^2_+)$ depending only
on $\ze^2_+$ {\it up to second derivatives}. Furthermore, we may
write $L_{\ze^2_+}(\ze^2_+,\ze^4)$ as a {\it linear operator}
$L_{\ze^2_+}$ applied to $(\ze^2_+,\ze^4)$, where the
coefficients of $L_{\ze^2_+}$ depend on $\ze^2_+$ only {\it up
to first derivatives}.

The operator $L_{\ze^2_+}$ is essentially $\De_C(\d_++\d^*)$,
and so is a {\it linear third-order elliptic operator}. Note
that we allow $0\le n\le m+1$ in \eq{co4eq9}. The reason for
this is that the {\it coefficients} of $L_{\ze^2_+}$ depend on
$\ze^2_+,\na\ze^2_+$, and so lie in $L^p_{m+1}$ locally. So the
maximum regularity we can expect for $L_{\ze^2_+}(\xi^2_+,\xi^4)$
is $L^p_{m+1}$, forcing~$n\le m+1$.

The most obvious value for $n$ in \eq{co4eq9} is $n=m-1$, as
then the domain of $L_{\ze^2_+}$ is the space $L^p_{m+2,\ga}$
containing $(\ze^2_+,\ze^4)$. However, our next lemma is an
{\it elliptic regularity result\/} for $L_{\ze^2_+}$ when $n=m$.
This is because we will use it to increase the regularity of
$(\ze^2_+,\ze^4)$ from $L^p_{m+2,\ga}$ to $L^p_{m+3,\ga}$ by
`bootstrapping'.

\begin{lem} There exists $A>0$ such that if\/ $(\xi^2_+,\xi^4)
\in L^p_{3,\ga}(\La^2_+T^*C\op\La^4T^*C)$ with\/
$L_{\ze^2_+}(\xi^2_+,\xi^4)\in L^p_{m,\ga}(\La^3T^*C)$ then
$(\xi^2_+,\xi^4)\in L^p_{m+3,\ga}(\La^2_+T^*C\op\La^4T^*C)$, and
\begin{equation}
\bnm{(\xi^2_+,\xi^4)}_{L^p_{m+3,\ga}}\le A\bigl(
\bnm{L_{\ze^2_+}(\xi^2_+,\xi^4)}_{L^p_{m,\ga}}
+\bnm{(\xi^2_+,\xi^4)}_{L^p_{0,\ga}}\bigl).
\label{co4eq11}
\end{equation}
\label{co4lem2}
\end{lem}

\begin{proof} If $L_{\ze^2_+}$ were a {\it smooth, asymptotically
cylindrical\/} elliptic operator the lemma would follow from
Lockhart and McOwen \cite[Theorem~3.7.2]{Lock}, \cite[Equation~(2.4)]{LoMc}.
We have to deal with the fact that the coefficients of $L_{\ze^2_+}$
depend on $\ze^2_+,\na\ze^2_+$, and so are only $L^p_{m+1}$ locally
rather than smooth. Unfortunately, Lockhart and McOwen's proof
\cite[p420]{LoMc} is not very informative, saying only that
\eq{co4eq11} is established `by standard parametric techniques'.

{\it Local\/} results of this form are proved by Morrey
\cite[Section~6.2]{Morr}. Supposing only that the coefficients
of the third-order elliptic operator $L_{\ze^2_+}$ are $C^m$,
Morrey \cite[Theorem~6.2.5]{Morr} implies that if $(\xi^2_+,\xi^4)$
is locally $L^p_3$ and $L_{\ze^2_+}(\xi^2_+,\xi^4)$ is locally
$L^p_m$ then $(\xi^2_+,\xi^4)$ is locally $L^p_{m+3}$. Also,
\cite[Theorem~6.2.6]{Morr} proves a local interior estimate of
the form \eq{co4eq11}, where $A>0$ depends on $m,p$, the
domains involved, $C^m$ bounds on the coefficients of
$L_{\ze^2_+}$, and a {\it modulus of continuity} for
their $m^{\rm th}$ derivatives.

Such a modulus of continuity is provided by a H\"older
$C^{0,\al}$ bound for the $m^{\rm th}$ derivatives, for
$\al\in(0,1)$. Thus, we can prove local estimates of
the form \eq{co4eq11} provided we have local $C^{m,\al}$
bounds for the coefficients of $L_{\ze^2_+}$, which
follow from local $C^{m+1,\al}$ bounds for $\ze^2_+$.
Now $\ze^2_+\in L^p_{m+2,\ga}(\La^2_+T^*C)$ and $p>4$,
so the Sobolev Embedding Theorem shows that $L^p_{m+2}$
embeds in $C^{m+1,\al}$ for $\al=1-4/p$. Therefore we
do have $C^{m,\al}$ control on the coefficients of
$L_{\ze^2_+}$, and can show that $L_{\ze^2_+}$ is
asymptotic in a weighted H\"older $C^{m,\al}_\ga$
sense to a cylindrical operator. So \eq{co4eq11} holds
as in \cite{LoMc}, and the lemma is proved.
\end{proof}

Suppose now that $\ze^2_+\in Q^{-1}(0)$. Then for $m=l\ge 1$ we
have $(\ze^2_+,0)\in L^p_{m+2,\ga}\bigl(B_{\ep'}(\La^2_+T^*C)
\bigr)\t L^p_{m+2,\ga}(\La^4T^*C)$, with $P(\ze^2_+,0)\equiv 0$.
So \eq{co4eq8} gives
\begin{equation}
L_{\ze^2_+}(\ze^2_+,0)=-E(x,\ze^2_+,\na\ze^2_+,\na^2\ze^2_+).
\label{co4eq12}
\end{equation}
As $\ze^2_+\in L^p_{m+2,\ga}(\La^2_+T^*C)$ we see from
\eq{co4eq10} and the asymptotic behavior of $F$ that the right
hand side of \eq{co4eq12} lies in $L^p_{m,\ga}(\La^3T^*C)$. Lemma
\ref{co4lem2} with $(\xi^2_+,\xi^4)=(\ze^2_+,0)$ then shows that
$\ze^2_+\in L^p_{m+3,\ga}(\La^2_+T^*C)$. So we have increased the
regularity of $\ze^2_+$ by one derivative. By induction we have
$\ze^2_+\in L^p_{m+2,\ga}(\La^2_+T^*C)$ for $m=l,l+1,l+2,\ldots$,
and Proposition \ref{co4prop4} is proved.
\end{proof}

Now let $\M_C^\ga$ be the moduli space of asymptotically
cylindrical coassociative submanifolds in $M$ close to $C$,
and asymptotic to $L\t(R',\iy)$ with decay rate $\ga$, as
in Theorem \ref{co1thm}. Define a map $S\co Q^{-1}(0)\ra
\{$4--submanifolds of $M\}$ by $S\co \ze^2_+\ra\Th(\Ga_{\ze^2_+})$,
where $\Ga_{\ze^2_+}$ is the {\it graph\/} of $\ze^2_+$ in
$B_{\ep'}(\La^2_+T^*C)$ as above.

\begin{prop} $S$ is a homeomorphism from $Q^{-1}(0)$ to
a neighborhood of\/ $C$ in~$\M_C^\ga$. \label{co4prop5}
\end{prop}

\begin{proof} First we must show that $S$ maps $Q^{-1}(0)\ra
\M_C^\ga$. Let $\ze^2_+\in Q^{-1}(0)$, and set $\ti C=S(\ze^2_+)$.
Then $(\Th\ci\ze^2_+)^*(\vp)=Q(\ze^2_+)=0$, so $\vp\vert_{\ti C}=0$
as $\ti C$ is the image of the map $\Th\ci\ze^2_+\co C\ra M$. Hence
$\ti C$ is {\it coassociative}. Proposition \ref{co4prop4} and
Sobolev Embedding shows that $\ze^2_+$ is smooth, so $\ti C$ is
a {\it smooth submanifold\/} of~$M$.

To show $\ti C$ is {\it asymptotically cylindrical} with rate
$\ga$, note that $\Th\ci\ze^2_+\co C\ra\ti C$ is a diffeomorphism.
Define $\ti K'=\Th\ci\ze^2_+(K')$ and $\ti\Phi\co L\t(R',\iy)\ra\ti
C\sm\ti K'$ by $\ti\Phi=\Th\ci\ze^2_+\ci\Phi$. Then $\ti K'$
is compact, and $\ti\Phi$ is a diffeomorphism. Now $\ze^2_+$ is
a section of $\La^2_+T^*C\cong\nu_C$, so $\Phi^*(\ze^2_+)$ is
a section of $\Phi^*(\nu_C)$ over $L\t(R',\iy)$. Pulling back
by the isomorphism $\xi$ gives a section $\xi^*\ci\Phi^*(\ze^2_+)$
of the vector bundle $\nu_L\t(R',\iy)$ over~$L\t(R',\iy)$.

Now the normal vector field $v$ on $L\t(R',\iy)$ making
\eq{co2eq3} commute is also a section of $\nu_L\t(R',\iy)$.
Define $\ti v=v+\xi^*\ci\Phi^*(\ze^2_+)$. Then $\ti v$ is a
section of $L\t(R',\iy)$, and the definition of $\Xi$ in
\eq{co4eq2} and part (iii) of the definition of $\Th$
after \eq{co4eq5} show that \eq{co2eq3} commutes with
$C,K',\Phi,v$ replaced by $\ti C,\ti K',\ti\Phi,\ti v$.
Therefore by Definition \ref{co2def6}, $\ti C$ is
asymptotically cylindrical with rate $\ga$ if $\bmd{\na^k
\ti v}=O(e^{\ga t})$ on $L\t(R',\iy)$ for all~$k\ge 0$.

As $C$ is asymptotically cylindrical with rate $\be<\ga$ we have
$\bmd{\na^kv}=O(e^{\be t})$ for all $k\ge 0$. Proposition
\ref{co4prop4} and Sobolev Embedding for weighted spaces, which
holds as in Lockhart \cite[Theorem~3.10]{Lock} and Bartnik
\cite[Theorem~1.2]{Bart}, then imply that
$\bmd{\na^k\ze^2_+}=O(e^{\ga\rho})$ on $C$ for all $k\ge 0$, where
$\rho$ is as in Definition \ref{co3def1}. Therefore
$\na^k(\xi^*\ci\Phi^*(\ze^2_+)) =O(e^{\ga t})$ on $L\t(R',\iy)$,
as $\Phi^*(\rho)\equiv t$, and $\xi,\Phi$ are asymptotically
cylindrical. As $\be<\ga$ this gives $\bmd{\na^k\ti v}=O(e^{\ga
t})$ for all $k\ge 0$, and $\ti C$ is asymptotically cylindrical
with weight $\ga$. Hence $S$ maps~$Q^{-1}(0)\ra\M_C^\ga$.

Next we reverse the argument. Suppose $\ti C$ is close to
$C$ in $\M_C^\ga$. As $\ti C,C$ are $C^1$ close there exists
a unique smooth section $\ze^2_+$ of $B_{\ep'}(\La^2_+T^*C)$
such that $\Th\ci\ze^2_+\co C\ra\ti C$ is a diffeomorphism.
Since $\vp\vert_{\ti C}\equiv 0$ we have $Q(\ze^2_+)=
(\Th\ci\ze^2_+)^*(\vp)=0$. Let $\ti C$ have data $\ti K'$,
$\ti\Phi$, $\ti v$ as in Definition \ref{co2def6}. The
argument above shows that $\ti v=v+\xi^*\ci\Phi^*(\ze^2_+)$.
But $\bmd{\na^kv}=O(e^{\ga t})$ and $\bmd{\na^k\ti v}=O(
e^{\ga t})$ for all $k\ge 0$. Subtracting implies that
$\bmd{\na^k\ze^2_+}=O(e^{\ga\rho})$ on $C$ for all~$k\ge 0$.

We need to show that $\ze^2_+\in L^p_{l+2,\ga}\bigl(B_{\ep'}
(\La^2_+T^*C)\bigr)$. The estimate $\bmd{\na^k\ze^2_+}=O(
e^{\ga\rho})$ for $k\ge 0$ does not prove this, however it
does imply that $\ze^2_+$ lies in $L^p_{l+2,\ga'}\bigl(B_{\ep'}
(\La^2_+T^*C)\bigr)$ for any $\ga'>\ga$. Hence $\ze^2_+$
lies in $Q'{}^{-1}(0)$, where $Q'$ is $Q$ defined with
rate~$\ga'$.

Choose $\ga'>\ga$ such that $[\ga,\ga']\cap\D_{(\d_++\d^*)_0}
=\emptyset$. Then $Q^{-1}(0)\subseteq Q'{}^{-1}(0)$, and from
above $Q^{-1}(0)$, $Q'{}^{-1}(0)$ are smooth, finite-dimensional
and locally isomorphic to $\Ker\bigl((\d_++\d^*)^p_{l+2,\ga}\bigr)$,
$\Ker\bigl((\d_++\d^*)^p_{l+2,\ga'}\bigr)$ respectively. But
as $[\ga,\ga']\cap\D_{(\d_++\d^*)_0}=\emptyset$ these kernels
are equal, so $Q^{-1}(0)$ and $Q'{}^{-1}(0)$ coincide near 0.
Thus if $\ti C$ is close enough to $C$ in $\M_C^\ga$ then
$\ze^2_+$ lies in $Q^{-1}(0)$, as we want.

To show $S$ is a {\it homeomorphism} requires us to
specify the topology on $\M_C^\ga$, which we have not done.
The natural way to define a topology on a space of submanifolds
is to identify submanifolds $\ti C$ near $C$ with sections
of the normal bundle $\nu_C$ to $C$, and induce the topology
from some Banach norm on a space of sections of $\nu_C$. In
our case, this just means that the topology on $\M_C^\ga$
is induced from some choice of Banach norm on sections
$\ze^2_+$ of $\La^2_+T^*C$, say the $C^1_\ga$ topology.

But since $Q^{-1}(0)$ with its $L^p_{l+2,\ga}$ topology
is locally homeomorphic to the {\it finite-dimensional\/}
vector space $\Ker\bigl((\d_++\d^*)^p_{l+2,\ga}\bigr)$,
{\it any} choice of topology on sections $\ze^2_+$ gives
the same topology on $\M_C^\ga$, as all Banach norms give
the same topology on a finite-dimensional space. So $S$
is a local homeomorphism.
\end{proof}

As from above $Q^{-1}(0)$ is smooth, finite-dimensional and
locally isomorphic to $\Ker\bigl((\d_++\d^*)^p_{l+2,\ga}\bigr)$,
Theorem \ref{co3thm2} and Proposition \ref{co4prop5} now prove
Theorem~\ref{co1thm}.

\begin{rem} In a forthcoming paper, we study the deformation space of
asymptotically cylindrical coassociative submanifolds with moving
boundary. Similarly to the setting here, we work with an
asymptotically cylindrical $G_2$--manifold $M$ with a Calabi--Yau
boundary $X$ at infinity. We also assume that the boundary of the
coassociative submanifold at infinity is a special Lagrangian
submanifold of $X$ and is allowed to move. The details of this
construction will appear in~\cite{Salu2}.
\end{rem}

\textbf{Acknowledgements}\qua Part of this work was done when the
second author was visiting Oxford University during Fall 2003.
Many thanks to the AWM for their grant support which made this
visit possible. We also would like to thank Conan Leung for
bringing this problem to our attention. The first author would
like to thank Jason Lotay for useful conversations.

\end{document}